# FAST LEARNING RATES IN STATISTICAL INFERENCE THROUGH AGGREGATION

By Jean-Yves Audibert

*Université Paris Est*

We develop minimax optimal risk bounds for the general learning task consisting in predicting as well as the best function in a reference set $\mathcal{G}$ up to the smallest possible additive term, called the convergence rate. When the reference set is finite and when $n$ denotes the size of the training data, we provide minimax convergence rates of the form $C(\frac{\log|\mathcal{G}|}{n})^v$ with tight evaluation of the positive constant $C$ and with exact $0 < v \leq 1$, the latter value depending on the convexity of the loss function and on the level of noise in the output distribution.

The risk upper bounds are based on a sequential randomized algorithm, which at each step concentrates on functions having both low risk and low variance with respect to the previous step prediction function. Our analysis puts forward the links between the probabilistic and worst-case viewpoints, and allows to obtain risk bounds unachievable with the standard statistical learning approach. One of the key ideas of this work is to use probabilistic inequalities with respect to appropriate (Gibbs) distributions on the prediction function space instead of using them with respect to the distribution generating the data.

The risk lower bounds are based on refinements of the Assouad lemma taking particularly into account the properties of the loss function. Our key example to illustrate the upper and lower bounds is to consider the $L_q$-regression setting for which an exhaustive analysis of the convergence rates is given while $q$ ranges in $[1;+\infty[$.

**1. Introduction.** We are given a family $\mathcal{G}$ of functions and we want to learn from data a function that predicts as well as the best function in $\mathcal{G}$ up to some additive term called the convergence rate. Even when the set $\mathcal{G}$ is finite, this learning task is crucial since:

- any continuous set of prediction functions can be viewed through its covering nets with respect to (w.r.t.) appropriate (pseudo-)distances and these nets are generally finite;
- one way of doing model selection among a finite family of submodels is to cut the training set into two parts, use the first part to learn the best prediction function of each submodel and use the second part to learn a prediction function which performs as well as the best of the prediction functions learned on the first part of the training set.

From this last item, our learning task for finite $\mathcal{G}$ is often referred to as model selection aggregation. It has two well-known variants. Instead of looking for a function predicting as well as the best in $\mathcal{G}$, these variants want to perform as well as the best convex combination of functions in $\mathcal{G}$ or as well as the best linear combination of functions in $\mathcal{G}$. These three aggregation tasks are linked in several ways (see [45] and references within).

Nevertheless, among these learning tasks, model selection aggregation has rare properties. First, in general an algorithm picking functions in the set $\mathcal{G}$ is not optimal (see, e.g., [9], Theorem 2, [40], Theorem 3, [21], page 14).

This means that the estimator has to look at an enlarged set of prediction functions. Second, in the statistical community, the only known optimal algorithms are all based on a Cesaro mean of Bayesian estimators (also referred to as progressive mixture rule). Third, the proof of their optimality is not achieved by the most prominent tool in statistical learning theory: bounds on the supremum of empirical processes (see [48], and refined works as [13, 17, 37, 42] and references within).

The idea of the proof, which comes back to Barron [11], is based on a chain rule and appeared to be successful for least square and entropy losses [12, 19, 20, 21, 53] and for general loss in [34].

In online prediction with expert advice setting, without any probabilistic assumption on the generation of the data, appropriate weighting methods have been shown to behave as well as the best expert up to a minimax-optimal additive remainder term (see [26, 43] and references within). In this worst-case context, amazingly sharp constants have been found (see in particular [24, 25, 33, 54]). These results are expressed in cumulative loss and can be transposed to model selection aggregation to the extent that the expected risk of the randomized procedure based on sequential predictions is proportional to the expectation of the cumulative loss of the sequential procedure (see Lemma 4.3 for precise statement).

This work presents a sequential algorithm, which iteratively updates a prior distribution put on the set of prediction functions. Contrary to previously mentioned works, these updates take into account the variance of the task. As a consequence, posterior distributions concentrate on simultaneously low risk functions and functions close to the previously drawn prediction function. This conservative law is not surprising in view of previous works on high-dimensional statistical tasks, such as wavelet thresholding, shrinkage procedures, iterative compression schemes [5] and iterative feature selection [1].

The paper is organized as follows. Section 2 introduces the notation and the existing algorithms. Section 3 proposes a unifying setting to combine worst-case analysis tight results and probabilistic tools. It details our sequentially randomized estimator and gives a sharp expected risk bound. In Sections 4 and 5, we show how to apply our main result under assumptions coming respectively from sequential prediction and model selection aggregation. While all this work concentrates on stating results when the data are independent and identically distributed, Section 4.2 shows that the argument underlying the main theorem can be applied for sequential predictions in which no probabilistic assumption is made and in which the data points come one by one (i.e., not in a batch manner). Section 6 contains algorithms that satisfy sharp standard-style generalization error bounds. To the author's knowledge, these bounds are not achievable with a classical statistical learning approach based on supremum of empirical processes. Here the main trick is to use probabilistic inequalities w.r.t. appropriate distributions on the prediction function space instead of using them w.r.t. the distribution generating the data. Section 7 presents an improved bound for $L_q$-regression ($q > 1$) when the noise has just a bounded moment of order $s \geq q$. This last assumption is much weaker than the traditional exponential moment assumption. Section 8 refines Assouad's lemma in order to obtain sharp constants and to take into account the properties of the loss function of the learning task. We illustrate our results by providing lower bounds matching the upper bounds obtained in the previous sections and by improving significantly the constants in lower bounds concerning Vapnik–Cervonenkis classes in classification. Section 9 summarizes the contributions of this work and lists some related open problems.

## 2. Notation and existing algorithms.

We assume that we observe $n$ pairs $Z_1 = (X_1, Y_1), \ldots, Z_n = (X_n, Y_n)$ of input–output and that each pair has been independently drawn from the same unknown distribution denoted $P$. The input and output space are denoted respectively $\mathcal{X}$ and $\mathcal{Y}$, so that $P$ is a probability distribution on the product space $\mathcal{Z} \triangleq \mathcal{X} \times \mathcal{Y}$. The target of a learning algorithm is to predict the output $Y$ associated with an input $X$ for pairs $(X, Y)$ drawn from the distribution $P$. In this work, $Z_{n+1}$ will denote a random variable independent of the training set $Z_1^n \triangleq (Z_1, \ldots, Z_n)$ and with the same distribution $P$. The quality of a prediction function $g : \mathcal{X} \to \mathcal{Y}$ is measured by the *risk* (also called expected loss or regret):

$$R(g) \triangleq \mathbb{E}_{Z \sim P} \, L(Z, g),$$

where $L(Z, g)$ assesses the loss of considering the prediction function $g$ on the data $Z \in \mathcal{Z}$. The symbol $\triangleq$ is used to underline that the equality is a definition. When there is no ambiguity on the distribution that a random variable has, the expectation w.r.t. this distribution will simply be written by indexing the expectation sign $\mathbb{E}$ by the random variable. For instance, we can write $R(g) \triangleq \mathbb{E}_Z L(Z, g)$. More generally, when they are multiple sources of randomness, $\mathbb{E}_Z$ means that we take the expectation with respect to the conditional distribution of $Z$ knowing all other sources of randomness.

We use $L(Z, g)$ rather than $L[Y, g(X)]$ to underline that our results are not restricted to nonregularized losses, where we call nonregularized loss a loss that can be written as $\ell[Y, g(X)]$ for some function $\ell : \mathcal{Y} \times \mathcal{Y} \to \mathbb{R}$.

For any $i \in \{0, \ldots, n\}$, the *cumulative loss* suffered by the prediction function $g$ on the first $i$ pairs of input–output, denoted $Z_1^i$ for short, is

$$\Sigma_i(g) \triangleq \sum_{j=1}^{i} L(Z_j, g),$$

where by convention we take $\Sigma_0$ identically equal to zero. The symbol $\equiv$ is used to underline when a function is identical to a constant (e.g., $\Sigma_0 \equiv 0$). With slight abuse, a symbol denoting a constant function may be used to denote the value of this function.

We assume that the set, denoted $\bar{\mathcal{G}}$, of all prediction functions has been equipped with a $\sigma$-algebra. Let $\mathcal{D}$ be the set of all probability distributions on $\bar{\mathcal{G}}$. By definition, a randomized algorithm produces a prediction function drawn according to a probability in $\mathcal{D}$. Let $\mathcal{P}$ be a set of probability distributions on $\mathcal{Z}$ in which we assume that the true unknown distribution generating the data lies. The learning task is essentially described by the 3-tuple $(\mathcal{G}, L, P)$ since we look for a possibly randomized estimator (or algorithm) $\hat{g}$ such that

$$\sup_{P \in \mathcal{P}} \left\{ \mathbb{E}_{Z_1^n} R(\hat{g}_{Z_1^n}) - \min_{g \in \mathcal{G}} R(g) \right\}$$

is minimized, where we recall that $R(g) \triangleq \mathbb{E}_{Z \sim P} L(Z, g)$. To shorten notation, when no confusion can arise, the dependence of $\hat{g}_{Z_1^n}$ w.r.t. the training sample $Z_1^n$ will be dropped and we will simply write $\hat{g}$. This means that we use the same symbol for both the algorithm and the prediction function produced by the algorithm on a training sample.

We implicitly assume that the quantities we manipulate are measurable; in particular, we assume that a prediction function is a measurable function from $\mathcal{X}$ to $\mathcal{Y}$, the mapping $(x, y, g) \mapsto L[(x, y), g]$ is measurable, the estimators considered in our lower bounds are measurable, . . . .

The $n$-fold product of a distribution $\mu$, which is the distribution of a vector consisting of $n$ i.i.d. realizations of $\mu$, is denoted $\mu^{\otimes n}$. For instance, the distribution of $(Z_1, \ldots, Z_n)$ is $P^{\otimes n}$.

The symbol $C$ will denote some positive constant whose value may differ from line to line. The set of nonnegative real numbers is denoted $\mathbb{R}_+ = [0; +\infty[$. We define $\lfloor x \rfloor$ as the largest integer $k$ such that $k \leq x$. To shorten notation, any finite sequence $a_1, \ldots, a_n$ will occasionally be denoted $a_1^n$. For instance, the training set is $Z_1^n$.

To handle possibly continuous set $\mathcal{G}$, we consider that $\mathcal{G}$ is a measurable space and that we have some *prior distribution* $\pi$ on it. The set of probability distributions on $\mathcal{G}$ will be denoted $\mathcal{M}$. The *Kullback–Leibler divergence* between a distribution $\rho \in \mathcal{M}$ and the prior distribution $\pi$ is

$$K(\rho, \pi) \triangleq \begin{cases} \mathbb{E}_{g \sim \rho} \log\left(\dfrac{\rho}{\pi}(g)\right), & \text{if } \rho \ll \pi, \\ +\infty, & \text{otherwise,} \end{cases}$$

where $\frac{\rho}{\pi}$ denotes the density of $\rho$ w.r.t. $\pi$ when it exists (i.e., $\rho \ll \pi$). For any $\rho \in \mathcal{M}$, we have $K(\rho, \pi) \geq 0$ and when $\pi$ is the uniform distribution on a finite set $\mathcal{G}$, we also have $K(\rho, \pi) \leq \log |\mathcal{G}|$. The Kullback–Leibler divergence satisfies the duality formula (see, e.g., [22], page 160): for any real-valued measurable function $h$ defined on $\mathcal{G}$,

$$(2.1) \qquad \inf_{\rho \in \mathcal{M}} \{\mathbb{E}_{g \sim \rho} h(g) + K(\rho, \pi)\} = -\log \mathbb{E}_{g \sim \pi} e^{-h(g)},$$

and that the infimum is reached for the *Gibbs distribution*

$$(2.2) \qquad \pi_{-h}(dg) \triangleq \frac{e^{-h(g)}}{\mathbb{E}_{g' \sim \pi} e^{-h(g')}} \cdot \pi(dg).$$

Intuitively, the Gibbs distribution $\pi_{-h}$ concentrates on prediction functions $g$ that are close to minimizing the function $h : \mathcal{G} \to \mathbb{R}$.

For any $\rho \in \mathcal{M}$, $\mathbb{E}_{g \sim \rho} g : x \mapsto \mathbb{E}_{g \sim \rho} g(x) = \int g(x) \rho(dg)$ is called a mixture of prediction functions. When $\mathcal{G}$ is finite, a mixture is simply a convex combination. Throughout this work, whenever we consider mixtures of prediction functions, we implicitly assume that $\mathbb{E}_{g \sim \rho} g(x)$ belongs to $\mathcal{Y}$ for any $x$ so that the mixture is a prediction function. This is typically the case when $\mathcal{Y}$ is an interval of $\mathbb{R}$.

We will say that the loss function is convex when the function $g \mapsto L(z, g)$ is convex for any $z \in \mathcal{Z}$, equivalently $L(z, \mathbb{E}_{g \sim \rho} g) \leq \mathbb{E}_{g \sim \rho} L(z, g)$ for any $\rho \in \mathcal{M}$ and $z \in \mathcal{Z}$. In this work, we do not assume the loss function to be convex except when it is explicitly mentioned.

The algorithm used to prove optimal convergence rates for several different losses (see, e.g., [12, 16, 19, 20, 21, 34, 53]) is the following:

ALGORITHM A. Let $\lambda > 0$. Predict according to $\frac{1}{n+1} \sum_{i=0}^{n} \mathbb{E}_{g \sim \pi_{-\lambda \Sigma_i}} g$, where we recall that $\Sigma_i$ maps a function $g \in \mathcal{G}$ to its cumulative loss up to time $i$.

In other words, for a new input $x$, the prediction of the output given by Algorithm A is $\frac{1}{n+1}\sum_{i=0}^n \int g(x)e^{-\lambda\Sigma_i(g)}\pi(dg)/\int e^{-\lambda\Sigma_i(g)}\pi(dg)$. Algorithm A has also been used with the classification loss. For this nonconvex loss, it has the same properties as the empirical risk minimizer on $\mathcal{G}$ [38, 39]. To give the optimal convergence rate, the parameter $\lambda$ and the distribution $\pi$ should be appropriately chosen. When $\mathcal{G}$ is finite, the estimator belongs to the convex hull of the set $\mathcal{G}$.

From Vovk, Haussler, Kivinen and Warmuth works [33, 51, 52] and the link between cumulative loss in online setting and expected risk in the batch setting (see Lemma 4.3), an "optimal" algorithm is:

ALGORITHM B.  Let $\lambda > 0$. For any $i \in \{0,\dots,n\}$, let $\hat{h}_i$ be a prediction function such that

$$\forall z \in \mathcal{Z} \qquad L(z,\hat{h}_i) \leq -\frac{1}{\lambda}\log \mathbb{E}_{g\sim\pi_{-\lambda\Sigma_i}}e^{-\lambda L(z,g)}.$$

If one of the $\hat{h}_i$ does not exist, the algorithm is said to fail. Otherwise it predicts according to $\frac{1}{n+1}\sum_{i=0}^n \hat{h}_i$.

In particular, for appropriate $\lambda > 0$, this algorithm does not fail when the loss function is the square loss (i.e., $L(z,g) = [y - g(x)]^2$) and when the output space is bounded. Algorithm B is based on the same Gibbs distribution $\pi_{-\lambda\Sigma_i}$ as Algorithm A. Besides, in [33], Example 3.13, it is shown that Algorithm A is not in general a particular case of Algorithm B, and that Algorithm B will not generally produce a prediction function in the convex hull of $\mathcal{G}$, unlike Algorithm A. In Sections 4 and 5, we will see how both algorithms are connected to the SeqRand algorithm presented in the next section.

**3. The algorithm and its generalization error bound.**  The aim of this section is to build an algorithm with the best possible minimax convergence rate. The algorithm relies on the following central condition for which we recall that $\mathcal{G}$ is a subset of the set $\bar{\mathcal{G}}$ of all prediction functions and that $\mathcal{M}$ and $\mathcal{D}$ are the sets of all probability distributions on respectively $\mathcal{G}$ and $\bar{\mathcal{G}}$.

For any $\lambda > 0$, let $\delta_\lambda$ be a real-valued function defined on $\mathcal{Z} \times \mathcal{G} \times \bar{\mathcal{G}}$ that satisfies the following inequality, which will be referred to as the *variance inequality*:

$$\forall \rho \in \mathcal{M} \ \exists \hat{\pi}(\rho) \in \mathcal{D}$$

$$\sup_{P\in\mathcal{P}}\{\mathbb{E}_{Z\sim P}\mathbb{E}_{g'\sim\hat{\pi}(\rho)}\log\mathbb{E}_{g\sim\rho}e^{\lambda[L(Z,g')-L(Z,g)-\delta_\lambda(Z,g,g')]}\} \leq 0.$$

The variance inequality is our probabilistic version of the generic algorithm condition in the online prediction setting (see [51], proof of Theorem

---

Input: $\lambda > 0$ and $\pi$ a distribution on the set $\mathcal{G}$.

1. Define $\hat{\rho}_0 \triangleq \hat{\pi}(\pi)$ in the sense of the variance inequality and draw a function $\hat{g}_0$ according to this distribution. Let $S_0(g) = 0$ for any $g \in \mathcal{G}$.
2. For any $i \in \{1, \ldots, n\}$, iteratively define

   (3.1) $\qquad S_i(g) \triangleq S_{i-1}(g) + L(Z_i, g) + \delta_\lambda(Z_i, g, \hat{g}_{i-1}) \qquad$ for any $g \in \mathcal{G}$.

   and

   $\qquad \hat{\rho}_i \triangleq \hat{\pi}(\pi_{-\lambda S_i}) \qquad$ in the sense of the variance inequality

   and draw a function $\hat{g}_i$ according to the distribution $\hat{\rho}_i$.
3. Predict with a function drawn according to the uniform distribution on the finite set $\{\hat{g}_0, \ldots, \hat{g}_n\}$.

   Conditionally to the training set, the distribution of the output prediction function will be denoted $\hat{\mu}$.

---

Fig. 1. *The SeqRand algorithm.*

1, or more explicitly in [33], page 11), in which we added the variance function $\delta_\lambda$. Our results will be all the sharper as this variance function is small. To make the variance inequality more readable, let us say for the moment that:

- Without any assumption on $\mathcal{P}$, for several usual "strongly" convex loss functions, we may take $\delta_\lambda \equiv 0$ provided that $\lambda$ is a small enough constant (see Section 4).
- The variance inequality can be seen as a "small expectation" inequality. The usual viewpoint is to control the quantity $L(Z, g)$ by its expectation w.r.t. $Z$ and a variance term. Here, roughly, $L(Z, g)$ is mainly controlled by $L(Z, g')$, where $g'$ is appropriately chosen through the choice of $\hat{\pi}(\rho)$, plus the additive term $\delta_\lambda$. By definition this additive term does not depend on the particular probability distribution generating the data and leads to empirical compensation.
- In the examples we will be interested in throughout this work, $\hat{\pi}(\rho)$ will be equal either to $\rho$ or to a Dirac distribution on some function, which is *not necessarily in* $\mathcal{G}$.
- For any loss function $L$, any set $\mathcal{P}$ and any $\lambda > 0$, one may choose $\delta_\lambda(Z, g, g') = \frac{\lambda}{2}[L(Z, g) - L(Z, g')]^2$ (see Section 6).

Our results concern the sequentially randomized algorithm described in Figure 1, which for sake of shortness we will call the SeqRand algorithm.

REMARK 3.1. When $\delta_\lambda(Z, g, g')$ does not depend on $g$, we recover a more standard-style algorithm to the extent that we then have $\pi_{-\lambda S_i} = \pi_{-\lambda \Sigma_i}$. Precisely our algorithm becomes the randomized version of Algorithm A. When $\delta_\lambda(Z, g, g')$ depends on $g$, the posterior distributions tend to concentrate on functions having small risk and small variance term. In Section 6, we will take $\delta_\lambda(Z, g, g') = \frac{\lambda}{2}[L(Z, g) - L(Z, g')]^2$. This choice implies a conservative mechanism: roughly, with high probability, among functions having low cumulative risk $\Sigma_i$, $\hat{g}_i$ will be chosen close to $\hat{g}_{i-1}$.

For any $i \in \{0, \ldots, n\}$, the quantities $S_i$, $\hat{\rho}_i$ and $\hat{g}_i$ depend on the training data only through $Z_1^i$, where we recall that $Z_1^i$ denotes $(Z_1, \ldots, Z_i)$. Besides they are also random to the extent that they depend on the draws of the functions $\hat{g}_0, \ldots, \hat{g}_{i-1}$.

The SeqRand algorithm produces a prediction function, which has three causes of randomness: the training data, the way $\hat{g}_i$ is obtained (step 2) and the uniform draw (step 3). For fixed $Z_1^i$ (i.e., conditional to $Z_1^i$), let $\Omega_i$ denote the joint distribution of $\hat{g}_0^i = (\hat{g}_0, \ldots, \hat{g}_i)$. The randomizing distribution $\hat{\mu}$ of the output prediction function by SeqRand is the distribution on $\bar{\mathcal{G}}$ corresponding to the last two causes of randomness. From the previous definitions, for any function $h : \bar{\mathcal{G}} \to \mathbb{R}$, we have $\mathbb{E}_{g \sim \hat{\mu}} h(g) = \mathbb{E}_{\hat{g}_0^n \sim \Omega_n} \frac{1}{n+1} \sum_{i=0}^n h(\hat{g}_i)$. Our main upper bound controls the expected risk $\mathbb{E}_{Z_1^n} \mathbb{E}_{g \sim \hat{\mu}} R(g)$ of the SeqRand procedure.

THEOREM 3.1. Let $\Delta_\lambda(g, g') \triangleq \mathbb{E}_{Z \sim P} \delta_\lambda(Z, g, g')$ for $g \in G$ and $g' \in \bar{\mathcal{G}}$, where we recall that $\delta_\lambda$ is a function satisfying the variance inequality. The expected risk of the SeqRand algorithm satisfies

$$
\begin{aligned}
(3.2) \quad & \mathbb{E}_{Z_1^n} \mathbb{E}_{g' \sim \hat{\mu}} R(g') \\
& \leq \min_{\rho \in \mathcal{M}} \left\{ \mathbb{E}_{g \sim \rho} R(g) + \mathbb{E}_{g \sim \rho} \mathbb{E}_{Z_1^n} \mathbb{E}_{g' \sim \hat{\mu}} \Delta_\lambda(g, g') + \frac{K(\rho, \pi)}{\lambda(n+1)} \right\}.
\end{aligned}
$$

In particular, when $\mathcal{G}$ is finite and when the loss function $L$ and the set $\mathcal{P}$ are such that $\delta_\lambda \equiv 0$, by taking $\pi$ uniform on $\mathcal{G}$, we get

$$
(3.3) \qquad \mathbb{E}_{Z_1^n} \mathbb{E}_{g \sim \hat{\mu}} R(g) \leq \min_{\mathcal{G}} R + \frac{\log |\mathcal{G}|}{\lambda(n+1)}.
$$

PROOF. Let $\mathcal{E}$ denote the expected risk of the SeqRand algorithm:

$$
\mathcal{E} \triangleq \mathbb{E}_{Z_1^n} \mathbb{E}_{g \sim \hat{\mu}} R(g) = \frac{1}{n+1} \sum_{i=0}^n \mathbb{E}_{Z_1^i} \mathbb{E}_{\hat{g}_0^i \sim \Omega_i} R(\hat{g}_i).
$$

We recall that $Z_{n+1}$ is a random variable independent of the training set $Z_1^n$ and with the same distribution $P$. Let $S_{n+1}$ be defined by (3.1) for $i = n+1$.

To shorten formulae, let $\hat{\pi}_i \triangleq \pi_{-\lambda S_i}$ so that by definition we have $\hat{\rho}_i = \hat{\pi}(\hat{\pi}_i)$. The variance inequality implies that

$$\mathbb{E}_{g' \sim \hat{\pi}(\rho)} R(g') \leq -\frac{1}{\lambda} \mathbb{E}_Z \mathbb{E}_{g' \sim \hat{\pi}(\rho)} \log \mathbb{E}_{g \sim \rho} e^{-\lambda[L(Z,g) + \delta_\lambda(Z,g,g')]}.$$

So for any $i \in \{0, \ldots, n\}$, for fixed $\hat{g}_0^{i-1} = (\hat{g}_0, \ldots, \hat{g}_{i-1})$ and fixed $Z_1^i$, we have

$$\mathbb{E}_{g' \sim \hat{\rho}_i} R(g') \leq -\frac{1}{\lambda} \mathbb{E}_{Z_{i+1}} \mathbb{E}_{g' \sim \hat{\rho}_i} \log \mathbb{E}_{g \sim \hat{\pi}_i} e^{-\lambda[L(Z_{i+1},g) + \delta_\lambda(Z_{i+1},g,g')]}.$$

Taking the expectations w.r.t. $(Z_1^i, \hat{g}_0^{i-1})$, we get

$$\mathbb{E}_{Z_1^i} \mathbb{E}_{\hat{g}_0^i} R(\hat{g}_i) = \mathbb{E}_{Z_1^i} \mathbb{E}_{\hat{g}_0^{i-1}} \mathbb{E}_{g' \sim \hat{\rho}_i} R(g')$$

$$\leq -\frac{1}{\lambda} \mathbb{E}_{Z_1^{i+1}} \mathbb{E}_{\hat{g}_0^i} \log \mathbb{E}_{g \sim \hat{\pi}_i} e^{-\lambda[L(Z_{i+1},g) + \delta_\lambda(Z_{i+1},g,\hat{g}_i)]}.$$

Consequently, by the chain rule (i.e., cancellation in the sum of logarithmic terms; [11]) and by intensive use of Fubini's theorem, we get

$$\mathcal{E} = \frac{1}{n+1} \sum_{i=0}^n \mathbb{E}_{Z_1^i} \mathbb{E}_{\hat{g}_0^i} R(\hat{g}_i)$$

$$\leq -\frac{1}{\lambda(n+1)} \sum_{i=0}^n \mathbb{E}_{Z_1^{i+1}} \mathbb{E}_{\hat{g}_0^i} \log \mathbb{E}_{g \sim \hat{\pi}_i} e^{-\lambda[L(Z_{i+1},g) + \delta_\lambda(Z_{i+1},g,\hat{g}_i)]}$$

$$= -\frac{1}{\lambda(n+1)} \mathbb{E}_{Z_1^{n+1}} \mathbb{E}_{\hat{g}_0^n} \sum_{i=0}^n \log \mathbb{E}_{g \sim \hat{\pi}_i} e^{-\lambda[L(Z_{i+1},g) + \delta_\lambda(Z_{i+1},g,\hat{g}_i)]}$$

$$= -\frac{1}{\lambda(n+1)} \mathbb{E}_{Z_1^{n+1}} \mathbb{E}_{\hat{g}_0^n} \sum_{i=0}^n \log \left( \frac{\mathbb{E}_{g \sim \pi} e^{-\lambda S_{i+1}(g)}}{\mathbb{E}_{g \sim \pi} e^{-\lambda S_i(g)}} \right)$$

$$= -\frac{1}{\lambda(n+1)} \mathbb{E}_{Z_1^{n+1}} \mathbb{E}_{\hat{g}_0^n} \log \left( \frac{\mathbb{E}_{g \sim \pi} e^{-\lambda S_{n+1}(g)}}{\mathbb{E}_{g \sim \pi} e^{-\lambda S_0(g)}} \right)$$

$$= -\frac{1}{\lambda(n+1)} \mathbb{E}_{Z_1^{n+1}} \mathbb{E}_{\hat{g}_0^n} \log \mathbb{E}_{g \sim \pi} e^{-\lambda S_{n+1}(g)}.$$

Now from the following lemma, we obtain

$$\mathcal{E} \leq -\frac{1}{\lambda(n+1)} \log \mathbb{E}_{g \sim \pi} e^{-\lambda \mathbb{E}_{Z_1^{n+1}} \mathbb{E}_{\hat{g}_0^n} S_{n+1}(g)}$$

$$= -\frac{1}{\lambda(n+1)} \log \mathbb{E}_{g \sim \pi} e^{-\lambda[(n+1)R(g) + \mathbb{E}_{Z_1^n} \mathbb{E}_{\hat{g}_0^n} \sum_{i=0}^n \Delta_\lambda(g,\hat{g}_i)]}$$

$$= \min_{\rho \in \mathcal{M}} \left\{ \mathbb{E}_{g \sim \rho} R(g) + \mathbb{E}_{g \sim \rho} \mathbb{E}_{Z_1^n} \mathbb{E}_{\hat{g}_0^n} \frac{\sum_{i=0}^n \Delta_\lambda(g,\hat{g}_i)}{n+1} + \frac{K(\rho,\pi)}{\lambda(n+1)} \right\}.$$

LEMMA 3.2. *Let $\mathcal{W}$ be a real-valued measurable function defined on a product space $\mathcal{A}_1 \times \mathcal{A}_2$ and let $\mu_1$ and $\mu_2$ be probability distributions on respectively $\mathcal{A}_1$ and $\mathcal{A}_2$ such that $\mathbb{E}_{a_1 \sim \mu_1} \log \mathbb{E}_{a_2 \sim \mu_2} e^{-\mathcal{W}(a_1, a_2)} < +\infty$. We have*

$$-\mathbb{E}_{a_1 \sim \mu_1} \log \mathbb{E}_{a_2 \sim \mu_2} e^{-\mathcal{W}(a_1, a_2)} \leq -\log \mathbb{E}_{a_2 \sim \mu_2} e^{-\mathbb{E}_{a_1 \sim \mu_1} \mathcal{W}(a_1, a_2)}.$$

PROOF. By using twice (2.1) and Fubini's theorem, we have

$$-\mathbb{E}_{a_1} \log \mathbb{E}_{a_2 \sim \mu_2} e^{-\mathcal{W}(a_1, a_2)} = \mathbb{E}_{a_1} \inf_{\rho} \{ \mathbb{E}_{a_2 \sim \rho} \mathcal{W}(a_1, a_2) + K(\rho, \mu_2) \}$$

$$\leq \inf_{\rho} \mathbb{E}_{a_1} \{ \mathbb{E}_{a_2 \sim \rho} \mathcal{W}(a_1, a_2) + K(\rho, \mu_2) \}$$

$$= -\log \mathbb{E}_{a_2 \sim \mu_2} e^{-\mathbb{E}_{a_1} \mathcal{W}(a_1, a_2)}. \qquad \square$$

Inequality (3.3) is a direct consequence of (3.2). $\square$

Theorem 3.1 bounds the expected risk of a randomized procedure, where the expectation is taken w.r.t. both the training set distribution and the randomizing distribution. From the following lemma, for convex loss functions, (3.3) implies

$$(3.4) \qquad \mathbb{E}_{Z_1^n} R(\mathbb{E}_{g \sim \hat{\mu}} g) \leq \min_{\mathcal{G}} R + \frac{\log |\mathcal{G}|}{\lambda(n+1)},$$

where we recall that $\hat{\mu}$ is the randomizing distribution of the SeqRand algorithm and $\lambda$ is a parameter whose typical value is the largest $\lambda > 0$ such that $\delta_\lambda \equiv 0$.

LEMMA 3.3. *For convex loss functions, the doubly expected risk of a randomized algorithm is greater than the expected risk of the deterministic version of the randomized algorithm; that is, if $\hat{\rho}$ denotes the randomizing distribution, we have*

$$\mathbb{E}_{Z_1^n} R(\mathbb{E}_{g \sim \hat{\rho}} g) \leq \mathbb{E}_{Z_1^n} \mathbb{E}_{g \sim \hat{\rho}} R(g).$$

PROOF. The result is a direct consequence of Jensen's inequality. $\square$

In [24], the authors rely on worst-case analysis to recover standard-style statistical results such as Vapnik's bounds [49]. Theorem 3.1 can be seen as a complement to this pioneering work. Inequality (3.4) is the model selection bound that is well known for least square regression and entropy loss, and that has been recently proved for general losses in [34].

Let us discuss the generalized form of the result. The right-hand side (r.h.s.) of (3.2) is a classical regularized risk, which appears naturally in the PAC-Bayesian approach (see, e.g., [7, 22, 56]). An advantage of stating the result this way is to be able to deal with uncountable infinite $\mathcal{G}$. Even when $\mathcal{G}$ is countable, this formulation has some benefit to the extent that for any measurable function $h : \mathcal{G} \to \mathbb{R}$, $\min_{\rho \in \mathcal{M}} \{ \mathbb{E}_{g \sim \rho} h(g) + K(\rho, \pi) \} \leq \min_{g \in \mathcal{G}} \{ h(g) + \log \pi^{-1}(g) \}$.

Our generalization error bounds depend on two quantities $\lambda$ and $\pi$ which are the parameters of our algorithm. Their choice depends on the precise setting. Nevertheless, when $\mathcal{G}$ is finite and with no particular structure a priori, a natural choice for $\pi$ is the uniform distribution on $\mathcal{G}$.

Once the distribution $\pi$ is fixed, an appropriate choice for the parameter $\lambda$ is the minimizer of the r.h.s. of (3.2). This minimizer is unknown by the statistician, and it is an open problem to adaptively choose $\lambda$ close to it.

**4. Link with sequential prediction.** This section aims at providing examples for which the variance inequality holds, at stating results coming from the online learning community in our batch setting (Section 4.1) and at providing new results for the sequential prediction setting in which no probabilistic assumption is made on the way the data are generated (Section 4.2).

4.1. *From online to batch.* In [33, 51, 52], the loss function is assumed to satisfy: there are positive numbers $\eta$ and $c$ such that

$$(4.1) \quad \begin{aligned} &\forall \rho \in \mathcal{M}, \; \exists g_\rho : \mathcal{X} \to \mathcal{Y}, \; \forall x \in \mathcal{X}, \; \forall y \in \mathcal{Y} \\ &L[(x,y), g_\rho] \leq -\frac{c}{\eta} \log \mathbb{E}_{g \sim \rho} e^{-\eta L[(x,y), g]}. \end{aligned}$$

REMARK 4.1. If $g \mapsto e^{-\eta L(z, g)}$ is concave, then (4.1) holds for $c = 1$ (and one may take $g_\rho = \mathbb{E}_{g \sim \rho} g$).

Assumption (4.1) implies that the variance inequality is satisfied both for $\lambda = \eta$ and $\delta_\lambda(Z, g, g') = (1 - 1/c) L(Z, g')$ and for $\lambda = \eta/c$ and $\delta_\lambda(Z, g, g') = (c - 1) L(Z, g)$, and we may take in both cases $\hat{\pi}(\rho)$ as the Dirac distribution at $g_\rho$. This leads to the *same* procedure that is described in the following straightforward corollary of Theorem 3.1.

COROLLARY 4.1. *Let $g_{\pi_{-\eta \Sigma_i}}$ be defined in the sense of (4.1) (for $\rho = \pi_{-\eta \Sigma_i}$). Consider the algorithm which predicts by drawing a function in $\{ g_{\pi_{-\eta \Sigma_0}}, \ldots, g_{\pi_{-\eta \Sigma_n}} \}$ according to the uniform distribution. Under assumption (4.1), its expected risk $\mathbb{E}_{Z_1^n} \frac{1}{n+1} \sum_{i=0}^{n} R(g_{\pi_{-\eta \Sigma_i}})$ is upper bounded by*

$$(4.2) \quad c \min_{\rho \in \mathcal{M}} \left\{ \mathbb{E}_{g \sim \rho} R(g) + \frac{K(\rho, \pi)}{\eta(n+1)} \right\}.$$

This result is not surprising in view of the following two results. The first one comes from worst-case analysis in sequential prediction.

THEOREM 4.2 ([33], Theorem 3.8). *Let $\mathcal{G}$ be countable. For any $g \in \mathcal{G}$, let $\Sigma_i(g) = \sum_{j=1}^{i} L(Z_j, g)$ (still) denote the cumulative loss up to time $i$ of the expert which always predicts according to function $g$. Under assumption (4.1), the cumulative loss on $Z_1^n$ of the strategy in which the prediction at time $i$ is done according to function $g_{\pi_{-\eta\Sigma_{i-1}}}$ in the sense of (4.1) (for $\rho = \pi_{-\eta\Sigma_{i-1}}$) is bounded by*

$$(4.3) \qquad \inf_{g \in \mathcal{G}} \left\{ c\Sigma_n(g) + \frac{c}{\eta} \log \pi^{-1}(g) \right\}.$$

The second result shows how the previous bound can be transposed into our model selection context by the following lemma.

LEMMA 4.3. *Let $\mathcal{A}$ be a learning algorithm which produces the prediction function $\mathcal{A}(Z_1^i)$ at time $i+1$, that is, from the data $Z_1^i = (Z_1, \ldots, Z_i)$. Let $\mathcal{L}$ be the randomized algorithm which produces a prediction function $\mathcal{L}(Z_1^n)$ drawn according to the uniform distribution on $\{\mathcal{A}(\varnothing), \mathcal{A}(Z_1), \ldots, \mathcal{A}(Z_1^n)\}$. The (doubly) expected risk of $\mathcal{L}$ is equal to $\frac{1}{n+1}$ times the expectation of the cumulative loss of $\mathcal{A}$ on the sequence $Z_1, \ldots, Z_{n+1}$.*

PROOF. By Fubini's theorem, we have

$$\mathbb{E}R[\mathcal{L}(Z_1^n)] = \frac{1}{n+1} \sum_{i=0}^{n} \mathbb{E}_{Z_1^n} R[\mathcal{A}(Z_1^i)]$$

$$= \frac{1}{n+1} \sum_{i=0}^{n} \mathbb{E}_{Z_1^{i+1}} L[Z_{i+1}, \mathcal{A}(Z_1^i)]$$

$$= \frac{1}{n+1} \mathbb{E}_{Z_1^{n+1}} \sum_{i=0}^{n} L[Z_{i+1}, \mathcal{A}(Z_1^i)]. \qquad \square$$

For any $\eta > 0$, let $c(\eta)$ denote the infimum of the $c$ for which (4.1) holds. Under weak assumptions, Vovk [52] proved that the infimum exists and studied the behavior of $c(\eta)$ and $a(\eta) = c(\eta)/\eta$, which are key quantities of (4.2) and (4.3). Under weak assumptions, and in particular in the examples given in Table 1, the optimal constants in (4.3) are $c(\eta)$ and $a(\eta)$ ([52], Theorem 1) and we have $c(\eta) \geq 1$, $\eta \mapsto c(\eta)$ nondecreasing and $\eta \mapsto a(\eta)$ nonincreasing. From these last properties, we understand the trade-off which occurs to choose the optimal $\eta$.

Table 1 specifies (4.2) in different well-known learning tasks. For instance, for bounded least square regression (i.e., when $|Y| \leq B$ for some $B > 0$),

TABLE 1
*Value of $c(\eta)$ for different loss functions*

| | Output space | Loss L(Z, g) | $c(\eta)$ |
|---|---|---|---|
| Entropy loss | $\mathcal{Y} = [0;1]$ | $Y\log(\frac{Y}{g(X)})$ | $c(\eta) = 1$ if $\eta \leq 1$ |
| [33], Example 4.3 | | $+(1-Y)\log(\frac{1-Y}{1-g(X)})$ | $c(\eta) = \infty$ if $\eta > 1$ |
| Absolute loss game | $\mathcal{Y} = [0;1]$ | $|Y - g(X)|$ | $\frac{\eta}{2\log[2/(1+e^{-\eta})]}$ |
| [33], Section 4.2 | | | $= 1 + \eta/4 + o(\eta)$ |
| Square loss | $\mathcal{Y} = [-B,B]$ | $[Y - g(X)]^2$ | $c(\eta) = 1$ if $\eta \leq 1/(2B^2)$ |
| [33], Example 4.4 | | | $c(\eta) = +\infty$ if $\eta > 1/(2B^2)$ |
| $L_q$-loss | $\mathcal{Y} = [-B,B]$ | $|Y - g(X)|^q$ | $c(\eta) = 1$ |
| (see Theorem 4.4) | | $q > 1$ | if $\eta \leq \frac{q-1}{qB^q}(1 \wedge 2^{2-q})$ |

Here $B$ denotes a positive real.

the generalization error of the algorithm described in Corollary 4.1 when $\eta = 1/(2B^2)$ is upper bounded by

$$(4.4) \qquad \min_{\rho \in \mathcal{M}} \left\{ \mathbb{E}_{g \sim \rho} R(g) + 2B^2 \frac{K(\rho, \pi)}{n+1} \right\}.$$

The constant appearing in front of the Kullback–Leibler divergence is much smaller than the ones obtained in unbounded regression setting even with Gaussian noise and bounded regression function (see [19, 34] and [22], page 87). The differences between these results partly come from the absence of boundedness assumptions on the output and from the weighted average used in the aforementioned works. Indeed the weighted average prediction function, that is, $\mathbb{E}_{g \sim \rho} g$, does not satisfy (4.1) for $c = 1$ and $\eta = 1/(2B^2)$ as was pointed out in [33], Example 3.13. Nevertheless, it satisfies (4.1) for $c = 1$ and $\eta \leq 1/(8B^2)$ (by using the concavity of $x \mapsto e^{-x^2}$ on $[-1/\sqrt{2}; 1/\sqrt{2}]$ and Remark 4.1), which leads to similar but weaker bound [see (4.2)].

*Case of the $L_q$-losses.* To deal with these losses, we need the following slight generalization of the result given in Appendix A of [35].

THEOREM 4.4. *Let $\mathcal{Y} = [a;b]$. We consider a nonregularized loss function, that is, a loss function such that $L(Z, g) = \ell[Y, g(X)]$ for any $Z = (X, Y) \in \mathcal{Z}$ and some function $\ell : \mathcal{Y} \times \mathcal{Y} \to \mathbb{R}$. For any $y \in \mathcal{Y}$, let $\ell_y$ be the function $[y' \mapsto \ell(y, y')]$. If for any $y \in \mathcal{Y}$:*

- *$\ell_y$ is continuous on $\mathcal{Y}$,*
- *$\ell_y$ decreases on $[a;y]$, increases on $[y;b]$ and $\ell_y(y) = 0$,*
- *$\ell_y$ is twice differentiable on the open set $(a;y) \cup (y;b)$,*

*then (4.1) is satisfied for $c = 1$ and*

$$(4.5) \qquad \eta \leq \inf_{a \leq y_1 < y < y_2 \leq b} \frac{\ell'_{y_1}(y)\ell''_{y_2}(y) - \ell''_{y_1}(y)\ell'_{y_2}(y)}{\ell'_{y_1}(y)[\ell'_{y_2}(y)]^2 - [\ell'_{y_1}(y)]^2\ell'_{y_2}(y)},$$

*where the infimum is taken w.r.t. $y_1, y$ and $y_2$.*

PROOF.   See Section 10.1.   □

REMARK 4.2.   This result simplifies the original one to the extent that $\ell_y$ does not need to be twice differentiable at point $y$ and the range of values for $y$ in the infimum is $(y_1; y_2)$ instead of $(a; b)$.

COROLLARY 4.5.   *For the $L_q$-loss, when $\mathcal{Y} = [-B; B]$ for some $B > 0$, condition (4.1) is satisfied for $c = 1$ and*

$$\eta \leq \frac{q-1}{qB^q}(1 \wedge 2^{2-q}).$$

PROOF.   We apply Theorem 4.4. By simple computations, the r.h.s. of (4.5) is

$$\inf_{-B \leq y_1 < y < y_2 \leq B} \frac{(q-1)(y_2 - y_1)}{q(y-y_1)(y_2-y)[(y-y_1)^{q-1} + (y_2-y)^{q-1}]}$$

$$= \frac{q-1}{q(2B)^q} \inf_{0 < t < 1} \frac{1}{t(1-t)[t^{q-1} + (1-t)^{q-1}]}.$$

For $1 < q \leq 2$, the infimum is reached for $t = 1/2$ and (4.5) can be written as $\eta \leq \frac{q-1}{qB^q}$. For $q \geq 2$, since the previous infimum is larger than $\inf_{0 < t < 1} \frac{1}{t(1-t)} = 4$, (4.5) is satisfied at least when $\eta \leq \frac{4(q-1)}{q(2B)^q}$.   □

4.2. *Sequential prediction.*   First note that using Corollary 4.5 and Theorem 4.2, we obtain a new result concerning sequential prediction for $L_q$-loss. Nevertheless this result is not due to our approach but to a refinement of the argument in [35], Appendix A. In this section, we will rather concentrate on giving results for sequential prediction coming from the arguments underlying Theorem 3.1.

In the online setting, the data points come one by one and there is no probabilistic assumption on the way they are generated. In this case, one should modify the definition of the variance function into: for any $\lambda > 0$, let $\delta_\lambda$ be a real-valued function defined on $\mathcal{Z} \times \mathcal{G} \times \bar{\mathcal{G}}$ that satisfies the following *online variance inequality*:

$$\forall \rho \in \mathcal{M}, \ \exists \hat{\pi}(\rho) \in \mathcal{D}, \ \forall z \in \mathcal{Z}$$

$$\mathbb{E}_{g' \sim \hat{\pi}(\rho)} \log \mathbb{E}_{g \sim \rho} e^{\lambda[L(z,g') - L(z,g) - \delta_\lambda(z,g,g')]} \leq 0.$$

---

Input: $\lambda > 0$ and $\pi$ a distribution on the set $\mathcal{G}$.

1. Define $\hat{\rho}_0 \triangleq \hat{\pi}(\pi)$ in the sense of the online variance inequality and draw a function $\hat{g}_0$ according to this distribution. For data $Z_1$, predict according to $\hat{g}_0$. Let $S_0(g) = 0$ for any $g \in \mathcal{G}$.

2. For any $i \in \{1, \ldots, n-1\}$, define

$$S_i(g) \triangleq S_{i-1}(g) + L(Z_i, g) + \delta_\lambda(Z_i, g, \hat{g}_{i-1}) \qquad \text{for any } g \in \mathcal{G},$$

and

$$\hat{\rho}_i \triangleq \hat{\pi}(\pi_{-\lambda S_i}) \qquad \text{in the sense of the online variance inequality}$$

and draw a function $\hat{g}_i$ according to the distribution $\hat{\rho}_i$. For data $Z_{i+1}$, predict according to $\hat{g}_i$.

---

Fig. 2. *The online SeqRand algorithm.*

The only difference with the variance inequality defined in Section 3 is the removal of the expectation with respect to $Z$. Naturally if $\delta_\lambda$ satisfies the online variance inequality, then it satisfies the variance inequality. The online version of the SeqRand algorithm is described in Figure 2. It satisfies the following theorem whose proof follows the same line as the one of Theorem 3.1.

Theorem 4.6. *The cumulative loss of the online SeqRand algorithm satisfies*

$$\sum_{i=1}^{n} \mathbb{E}_{\hat{g}_{i-1}} L(Z_i, \hat{g}_{i-1})$$

$$\leq \min_{\rho \in \mathcal{M}} \left\{ \mathbb{E}_{g \sim \rho} \sum_{i=1}^{n} L(Z_i, g) + \mathbb{E}_{g \sim \rho} \mathbb{E}_{\hat{g}_0^{n-1}} \sum_{i=1}^{n} \delta_\lambda(Z_i, g, \hat{g}_{i-1}) + \frac{K(\rho, \pi)}{\lambda} \right\}.$$

*In particular, when $\mathcal{G}$ is finite, by taking $\pi$ uniform on $\mathcal{G}$, we get*

$$\sum_{i=1}^{n} \mathbb{E}_{\hat{g}_{i-1}} L(Z_i, \hat{g}_{i-1})$$

$$\leq \min_{g \in \mathcal{G}} \left\{ \sum_{i=1}^{n} L(Z_i, g) + \mathbb{E}_{\hat{g}_0^{n-1}} \sum_{i=1}^{n} \delta_\lambda(Z_i, g, \hat{g}_{i-1}) + \frac{\log |\mathcal{G}|}{\lambda} \right\}.$$

Up to the online variance function $\delta_\lambda$, the online variance inequality is the generic algorithm condition of [33], page 11. So cases where $\delta_\lambda$ are equal to zero are already known. Now new results can be obtained by using that for any loss function $L$ and any $\lambda > 0$, the online variance inequality is satisfied for $\delta_\lambda(Z, g, g') = \frac{\lambda}{2}[L(Z, g) - L(Z, g')]^2$ (proof in Section 10.2). The associated distribution $\hat{\pi}(\rho)$ is then just $\rho$. In spirit, the result associated with these choices is similar to the ones obtained in [27], Section 4, to the extent that it gives a bound with second-order terms. Nevertheless, we do not know how to properly choose the parameter $\lambda$ whereas the aforementioned work solves this problem. More discussion on this topic can be found in [8], Section 4.2.

## 5. Model selection aggregation under Juditsky, Rigollet and Tsybakov assumptions [34].
The main result of [34] relies on the following assumption on the loss function $L$ and the set $\mathcal{P}$ of probability distributions on $\mathcal{Z}$ in which we assume that the true distribution lies. There exist $\lambda > 0$ and a real-valued function $\psi$ defined on $\mathcal{G} \times \mathcal{G}$ such that for any $P \in \mathcal{P}$

$$(5.1) \quad \begin{cases} \mathbb{E}_{Z \sim P} e^{\lambda[L(Z, g') - L(Z, g)]} \leq \psi(g', g), & \text{for any } g, g' \in \mathcal{G}, \\ \psi(g, g) = 1, & \text{for any } g \in \mathcal{G}, \\ \text{the function } [g \mapsto \psi(g', g)] \text{ is concave for any } g' \in \mathcal{G}. \end{cases}$$

Theorem 3.1 gives the following result.

COROLLARY 5.1. *Consider the algorithm which draws uniformly its prediction function in the set $\{\mathbb{E}_{g \sim \pi_{-\lambda \Sigma_0}} g, \ldots, \mathbb{E}_{g \sim \pi_{-\lambda \Sigma_n}} g\}$. Under assumption (5.1), its expected risk $\mathbb{E}_{Z_1^n} \frac{1}{n+1} \sum_{i=0}^n R(\mathbb{E}_{g \sim \pi_{-\lambda \Sigma_i}} g)$ is upper bounded by*

$$(5.2) \quad \min_{\rho \in \mathcal{M}} \left\{ \mathbb{E}_{g \sim \rho} R(g) + \frac{K(\rho, \pi)}{\lambda(n+1)} \right\}.$$

PROOF. We start by proving that the variance inequality holds with $\delta_\lambda \equiv 0$, and that we may take $\hat{\pi}(\rho)$ as the Dirac distribution at the function $\mathbb{E}_{g \sim \rho} g$. By using Jensen's inequality and Fubini's theorem, assumption (5.1) implies that

$$\mathbb{E}_{g' \sim \hat{\pi}(\rho)} \mathbb{E}_{Z \sim P} \log \mathbb{E}_{g \sim \rho} e^{\lambda[L(Z, g') - L(Z, g)]}$$

$$= \mathbb{E}_{Z \sim P} \log \mathbb{E}_{g \sim \rho} e^{\lambda[L(Z, \mathbb{E}_{g' \sim \rho} g') - L(Z, g)]}$$

$$\leq \log \mathbb{E}_{g \sim \rho} \mathbb{E}_{Z \sim P} e^{\lambda[L(Z, \mathbb{E}_{g' \sim \rho} g') - L(Z, g)]}$$

$$\leq \log \mathbb{E}_{g \sim \rho} \psi(\mathbb{E}_{g' \sim \rho} g', g)$$

$$\leq \log \psi(\mathbb{E}_{g' \sim \rho} g', \mathbb{E}_{g \sim \rho} g)$$

$$= 0,$$

so that we can apply Theorem 3.1. It remains to note that in this context the SeqRand algorithm is the one described in the corollary.  □

In this context, the SeqRand algorithm reduces to the randomized version of Algorithm A. From Lemma 3.3, for convex loss functions, (5.2) also holds for the risk of Algorithm A. Corollary 5.1 also shows that the risk bounds for Algorithm A proved in [34], Theorem 3.2, and the examples of [34], Section 4.2, hold with the same constants for the SeqRand algorithm (provided that the expected risk w.r.t. the training set distribution is replaced by the expected risk w.r.t. both training set and randomizing distributions).

On assumption (5.1) we should say that it does not a priori require the function $L$ to be convex. Nevertheless, any known relevant examples deal with "strongly" convex loss functions and we know that in general the assumption will not hold for the Support Vector Machine (or hinge loss) function and for the absolute loss function. Indeed, without further assumption, one cannot expect rates better than $1/\sqrt{n}$ for these loss functions (see Section 8.3).

By taking the appropriate variance function $\delta_\lambda(Z, g, g')$, it is possible to prove that the results in [34], Theorem 3.1, and [34], Section 4.1, hold for the SeqRand algorithm (provided that the expected risk w.r.t. the training set distribution is replaced by the expected risk w.r.t. both training set and randomizing distributions). The choice of $\delta_\lambda(Z, g, g')$, which for sake of shortness we do not specify, is in fact such that the resulting SeqRand algorithm is again the randomized version of Algorithm A.

**6. Standard-style statistical bounds.** This section proposes new results of a different kind. In the previous sections, under convexity assumptions, we were able to achieve fast rates. Here we have assumption neither on the loss function nor on the probability generating the data. Nevertheless we show that the SeqRand algorithm applied for $\delta_\lambda(Z, g, g') = \lambda[L(Z, g) - L(Z, g')]^2/2$ satisfies a sharp standard-style statistical bound.

This section contains two parts: the first one provides results in expectation (as in the preceding sections) whereas the second part provides deviation inequalities on the risk that require advances on the sequential prediction analysis.

6.1. *Bounds on the expected risk.*

6.1.1. *Bernstein's type bound.*

THEOREM 6.1. *Let $V(g, g') = \mathbb{E}_Z\{[L(Z, g) - L(Z, g')]^2\}$. Consider the SeqRand algorithm applied with $\delta_\lambda(Z, g, g') = \lambda[L(Z, g) - L(Z, g')]^2/2$ and $\hat{\pi}(\rho) = \rho$. Its expected risk $\mathbb{E}_{Z_1^n}\mathbb{E}_{g \sim \hat{\mu}}R(g)$, where we recall that $\hat{\mu}$ denotes the randomizing distribution, satisfies*

$$\mathbb{E}_{Z_1^n}\mathbb{E}_{g' \sim \hat{\mu}}R(g')$$

(6.1)
$$\leq \min_{\rho \in \mathcal{M}} \left\{ \mathbb{E}_{g \sim \rho} R(g) + \frac{\lambda}{2} \mathbb{E}_{g \sim \rho} \mathbb{E}_{Z_1^n} \mathbb{E}_{g' \sim \hat{\mu}} V(g, g') + \frac{K(\rho, \pi)}{\lambda(n+1)} \right\}.$$

PROOF. See Section 10.2.  □

To make (6.1) more explicit and to obtain a generalization error bound in which the randomizing distribution does not appear in the r.h.s. of the bound, the following corollary considers a widely used assumption relating the variance term to the excess risk (see Mammen and Tsybakov [41, 47], and also Polonik [44]). Precisely, from Theorem 6.1, we obtain:

COROLLARY 6.2. *If there exist* $0 \leq \gamma \leq 1$ *and a prediction function* $\tilde{g}$ *(not necessarily in* $\mathcal{G}$*) such that* $V(g, \tilde{g}) \leq c[R(g) - R(\tilde{g})]^\gamma$ *for any* $g \in \mathcal{G}$*, the expected risk* $\mathcal{E} = \mathbb{E}_{Z_1^n} \mathbb{E}_{g \sim \hat{\mu}} R(g)$ *of the SeqRand algorithm used in Theorem 6.1 satisfies:*

- *When* $\gamma = 1$,

$$\mathcal{E} - R(\tilde{g}) \leq \min_{\rho \in \mathcal{M}} \left\{ \frac{1 + c\lambda}{1 - c\lambda} [\mathbb{E}_{g \sim \rho} R(g) - R(\tilde{g})] + \frac{K(\rho, \pi)}{(1 - c\lambda)\lambda(n+1)} \right\}.$$

  *In particular, for* $\mathcal{G}$ *finite,* $\pi$ *the uniform distribution,* $\lambda = 1/(2c)$*, when* $\tilde{g}$ *belongs to* $\mathcal{G}$*, we get* $\mathcal{E} \leq \min_{g \in \mathcal{G}} R(g) + \frac{4c \log |\mathcal{G}|}{n+1}$.

- *When* $\gamma < 1$*, for any* $0 < \beta < 1$ *and for* $\tilde{R}(g) \triangleq R(g) - R(\tilde{g})$,

$$\mathcal{E} - R(\tilde{g}) \leq \left\{ \frac{1}{\beta} \min_{\rho \in \mathcal{M}} \left( \mathbb{E}_{g \sim \rho} [\tilde{R}(g) + c\lambda \tilde{R}^\gamma(g)] + \frac{K(\rho, \pi)}{\lambda(n+1)} \right) \right\}$$
$$\vee \left( \frac{c\lambda}{1 - \beta} \right)^{1/(1-\gamma)}.$$

PROOF. See Section 10.3.  □

To understand the sharpness of Theorem 6.1, we have to compare this result with the following one that comes from the traditional (PAC-Bayesian) statistical learning approach which relies on supremum of empirical processes. In the following theorem, we consider the estimator minimizing the uniform bound, that is, the estimator for which we have the smallest upper bound on its generalization error.

THEOREM 6.3. *We still use* $V(g, g') = \mathbb{E}_Z \{[L(Z, g) - L(Z, g')]^2\}$*. The generalization error of the algorithm which draws its prediction function* *according to the Gibbs distribution $\pi_{-\lambda \Sigma_n}$ satisfies*

$$\mathbb{E}_{Z_1^n} \mathbb{E}_{g' \sim \pi_{-\lambda \Sigma_n}} R(g')$$

$$(6.2) \quad \leq \min_{\rho \in \mathcal{M}} \Bigg\{ \mathbb{E}_{g \sim \rho} R(g) + \frac{K(\rho, \pi) + 1}{\lambda n} + \lambda \mathbb{E}_{g \sim \rho} \mathbb{E}_{Z_1^n} \mathbb{E}_{g' \sim \pi_{-\lambda \Sigma_n}} V(g, g')$$

$$+ \lambda \frac{1}{n} \sum_{i=1}^n \mathbb{E}_{g \sim \rho} \mathbb{E}_{Z_1^n} \mathbb{E}_{g' \sim \pi_{-\lambda \Sigma_n}} [L(Z_i, g) - L(Z_i, g')]^2 \Bigg\}.$$

*Let $\varphi$ be the positive convex increasing function defined as $\varphi(t) \triangleq \frac{e^t - 1 - t}{t^2}$ and $\varphi(0) = \frac{1}{2}$ by continuity. When $\sup_{z \in \mathcal{Z}, g \in \mathcal{G}, g' \in \mathcal{G}} |L(z, g') - L(z, g)| \leq B$, we also have*

$$\mathbb{E}_{Z_1^n} \mathbb{E}_{g' \sim \pi_{-\lambda \Sigma_n}} R(g')$$

$$(6.3) \quad \leq \min_{\rho \in \mathcal{M}} \Bigg\{ \mathbb{E}_{g \sim \rho} R(g) + \lambda \varphi(\lambda B) \mathbb{E}_{g \sim \rho} \mathbb{E}_{Z_1^n} \mathbb{E}_{g' \sim \pi_{-\lambda \Sigma_n}} V(g, g')$$

$$+ \frac{K(\rho, \pi) + 1}{\lambda n} \Bigg\}.$$

PROOF. See Section 10.4. □

As in Theorem 6.1, there is a variance term in which the randomizing distribution is involved. As in Corollary 6.2, one can convert (6.3) into a proper generalization error bound, that is, a nontrivial bound $\mathbb{E}_{Z_1^n} \mathbb{E}_{g \sim \pi_{-\lambda \Sigma_n}} R(g) \leq \mathcal{B}(n, \pi, \lambda)$ where the training data do not appear in $\mathcal{B}(n, \pi, \lambda)$.

By comparing (6.3) and (6.1), we see that the classical approach requires the quantity $\sup_{g \in \mathcal{G}, g' \in \mathcal{G}} |L(Z, g') - L(Z, g)|$ to be uniformly bounded and the unpleasing function $\varphi$ appears. In fact, using technical small expectations theorems (see, e.g., [4], Lemma 7.1), exponential moments conditions on the above quantity would be sufficient.

The symmetrization trick used to prove Theorem 6.1 is performed in the prediction functions space. We do not call on the second virtual training set currently used in statistical learning theory (see [49]). Nevertheless both symmetrization tricks end up to the same nice property: we need no boundedness assumption on the loss functions. In our setting, symmetrization on training data leads to an unwanted expectation and to a constant four times larger (see the two variance terms of (6.2) and the discussion in [5], Section 8.3.3).

In particular, deducing from Theorem 6.3 a corollary similar to Corollary 6.2 is only possible through (6.3) and provided that we have a boundedness assumption on $\sup_{z \in \mathcal{Z}, g \in \mathcal{G}, g' \in \mathcal{G}} |L(z, g') - L(z, g)|$. Indeed one cannot use (6.2) because of the last variance term in (6.2) (since $\Sigma_n$ depends on $Z_i$).

Our approach has nevertheless the following limit: the proof of Corollary 6.2 does not use a chaining argument. As a consequence, in the particular case when the model has polynomial entropies (see, e.g., [41]) and when the assumption in Corollary 6.2 holds for $\gamma < 1$ (and not for $\gamma = 1$), Corollary 6.2 does not give the minimax optimal convergence rate. Combining the better variance control presented here with the chaining argument is an open problem.

6.1.2. *Hoeffding's type bound.* Contrary to generalization error bounds coming from Bernstein's inequality, (6.1) does not require any boundedness assumption. For bounded losses, without any variance assumption (i.e., roughly when the assumption used in Corollary 6.2 does not hold for $\gamma > 0$), tighter results are obtained by using Hoeffding's inequality, that is: for any random variable $W$ satisfying $a \leq W \leq b$, then for any $\lambda > 0$

$$\mathbb{E}e^{\lambda(W - \mathbb{E}W)} \leq e^{\lambda^2(b-a)^2/8}.$$

THEOREM 6.4. *Assume that for any $z \in \mathcal{Z}$ and $g \in \mathcal{G}$, we have $a \leq L(z, g) \leq b$ for some reals $a, b$. Consider the SeqRand algorithm applied with $\delta_\lambda(Z, g, g') = \lambda(b-a)^2/8$ and $\hat{\pi}(\rho) = \rho$. Its expected risk $\mathbb{E}_{Z_1^n} \mathbb{E}_{g \sim \hat{\mu}} R(g)$, where we recall that $\hat{\mu}$ denotes the randomizing distribution, satisfies*

$$(6.4) \qquad \mathbb{E}_{Z_1^n} \mathbb{E}_{g \sim \hat{\mu}} R(g) \leq \min_{\rho \in \mathcal{M}} \left\{ \mathbb{E}_{g \sim \rho} R(g) + \frac{\lambda(b-a)^2}{8} + \frac{K(\rho, \pi)}{\lambda(n+1)} \right\}.$$

*In particular, when $\mathcal{G}$ is finite, by taking $\pi$ uniform on $\mathcal{G}$ and $\lambda = \sqrt{\frac{8 \log|\mathcal{G}|}{(b-a)^2(n+1)}}$, we get*

$$(6.5) \qquad \mathbb{E}_{Z_1^n} \mathbb{E}_{g \sim \hat{\mu}} R(g) - \min_{g \in \mathcal{G}} R(g) \leq (b-a)\sqrt{\frac{\log|\mathcal{G}|}{2(n+1)}}.$$

PROOF. From Hoeffding's inequality, we have

$$\mathbb{E}_{g' \sim \hat{\pi}(\rho)} \log \mathbb{E}_{g \sim \rho} e^{\lambda[L(Z, g') - L(Z, g)]} = \log \mathbb{E}_{g \sim \rho} e^{\lambda[\mathbb{E}_{g' \sim \hat{\pi}(\rho)} L(Z, g') - L(Z, g)]}$$

$$\leq \frac{\lambda^2(b-a)^2}{8},$$

hence the variance inequality holds for $\delta_\lambda \equiv \lambda(b-a)^2/8$ and $\hat{\pi}(\rho) = \rho$. The result directly follows from Theorem 3.1. $\square$

The standard point of view (see Appendix A.2) applies Hoeffding's inequality to the random variable $W = L(Z, g') - L(Z, g)$ for $g$ and $g'$ fixed and $Z$ drawn according to the probability generating the data. The previous theorem uses it on the random variable $W = L(Z, g') - \mathbb{E}_{g \sim \rho} L(Z, g)$ for fixed $Z$ and fixed probability distribution $\rho$ but for $g'$ drawn according to $\rho$. Here the gain is a multiplicative factor equal to 2 (see Appendix A.2).

6.2. *Deviation inequalities.* For the comparison between Theorem 6.1 and Theorem 6.3 to be fair, one should add that (6.3) and (6.2) come from deviation inequalities that are not exactly obtainable to the author's knowledge with the arguments developed here. Precisely, consider the following adaptation of Lemma 5 of [55].

LEMMA 6.5. *Let $\mathcal{A}$ be a learning algorithm which produces the prediction function $\mathcal{A}(Z_1^i)$ at time $i + 1$, that is, from the data $Z_1^i = (Z_1, \ldots, Z_i)$. Let $\mathcal{L}$ be the randomized algorithm which produces a prediction function $\mathcal{L}(Z_1^n)$ drawn according to the uniform distribution on $\{\mathcal{A}(\varnothing), \mathcal{A}(Z_1), \ldots, \mathcal{A}(Z_1^n)\}$. Assume that $\sup_{z, g, g'} |L(z, g) - L(z, g')| \leq B$ for some $B > 0$. Conditionally to $Z_1, \ldots, Z_{n+1}$, the expectation of the risk of $\mathcal{L}$ w.r.t. to the uniform draw is $\frac{1}{n+1} \sum_{i=0}^{n} R[\mathcal{A}(Z_1^i)]$ and satisfies: for any $\eta > 0$ and $\varepsilon > 0$, for any reference prediction function $\tilde{g}$, with probability at least $1 - \varepsilon$ w.r.t. the distribution of $Z_1, \ldots, Z_{n+1}$,*

$$
\begin{aligned}
(6.6) \qquad & \frac{1}{n+1} \sum_{i=0}^{n} R[\mathcal{A}(Z_1^i)] - R(\tilde{g}) \\
& \leq \frac{1}{n+1} \sum_{i=0}^{n} \{L[Z_{i+1}, \mathcal{A}(Z_1^i)] - L(Z_{i+1}, \tilde{g})\} \\
& \quad + \eta \varphi(\eta B) \frac{1}{n+1} \sum_{i=0}^{n} V[\mathcal{A}(Z_1^i), \tilde{g}] + \frac{\log(\varepsilon^{-1})}{\eta(n+1)},
\end{aligned}
$$

*where we still use $V(g, g') = \mathbb{E}_Z\{[L(Z, g) - L(Z, g')]^2\}$ for any prediction functions $g$ and $g'$ and $\varphi(t) \triangleq \frac{e^t - 1 - t}{t^2}$ for any $t > 0$.*

PROOF. See Section 10.5. □

We see that two variance terms appear. The first one comes from the worst-case analysis and is hidden in $\sum_{i=0}^{n} \{L[Z_{i+1}, \mathcal{A}(Z_1^i)] - L(Z_{i+1}, \tilde{g})\}$ and the second one comes from the concentration result (Lemma 10.1). The presence of this last variance term annihilates the benefits of our approach in which we were manipulating variance terms much smaller than the traditional Bernstein's variance term.

To illustrate this point, consider for instance least square regression with bounded outputs: from Theorem 4.2 and Table 1, the hidden variance term is null. In some situations, the second variance term $\frac{1}{n+1}\sum_{i=0}^{n}V[\mathcal{A}(Z_1^i),\tilde{g}]$ may behave like a positive constant; for instance, this occurs when $\mathcal{G}$ contains two very different functions having the optimal risk $\min_{g\in\mathcal{G}}R(g)$. By optimizing $\eta$, this will lead to a deviation inequality of order $n^{-1/2}$ even though from (4.4) the procedure has $n^{-1}$-convergence rate in expectation. In [9], Theorem 3, in a rather general learning setting, this deviation inequality of order $n^{-1/2}$ is proved to be optimal.

To conclude, for deviation inequalities, we cannot expect to do better than the standard-style approach since at some point we use a Bernstein's type bound w.r.t. the distribution generating the data. Besides procedures based on worst-case analysis seem to suffer higher fluctuations of the risk than necessary (see [9], discussion of Theorem 3).

REMARK 6.1.   Lemma 6.5 should be compared with Lemma 4.3. The latter deals with results in expectation while the former concerns deviation inequalities. Note that Lemma 6.5 requires the loss function to be bounded and makes a variance term appear.

**7. Application to $L_q$-regression for unbounded outputs.**   In this section, we consider the $L_q$-loss: $L(Z,g)=|Y-g(X)|^q$. As a warm-up exercise, we tackle the absolute loss setting (i.e., $q=1$). The following corollary holds without any assumption on the output (except naturally that if $E_Z|Y|<+\infty$ to ensure finite risk).

COROLLARY 7.1.   *Let $q=1$. Assume that $\sup_{g\in\mathcal{G}}E_Z\,g(X)^2\le b^2$ for some $b>0$. There exists an estimator $\hat{g}$ such that*

$$(7.1)\qquad\qquad \mathbb{E}R(\hat{g})-\min_{g\in\mathcal{G}}R(g)\le 2b\sqrt{\frac{2\log|\mathcal{G}|}{n+1}}.$$

PROOF.   Using $\mathbb{E}_Z\{[|Y-g(X)|-|Y-g'(X)|]^2\}\le 4b^2$ and Theorem 6.1, the algorithm considered in Theorem 6.1 satisfies $\mathbb{E}R(\hat{g})-\min_{g\in\mathcal{G}}R(g)\le 2\lambda b^2+\frac{\log|\mathcal{G}|}{\lambda(n+1)}$, which gives the desired result by taking $\lambda=\sqrt{\frac{\log|\mathcal{G}|}{2b^2(n+1)}}$.   □

Now we deal with the strongly convex loss functions (i.e., $q>1$). Using Theorem 3.1 jointly with the symmetrization idea developed in the previous section allows to obtain new convergence rates in heavy noise situation, that is, when the output is not constrained to have a bounded exponential moment.

COROLLARY 7.2. *Let $q > 1$. Assume that*

$$\begin{cases} \sup_{g \in \mathcal{G}, x \in \mathcal{X}} |g(x)| \leq b, & \text{for some } b > 0, \\ \mathbb{E}|Y|^s \leq A, & \text{for some } s \geq q \text{ and } A > 0, \\ \mathcal{G} \text{ finite.} \end{cases}$$

*Let $\pi$ be the uniform distribution on $\mathcal{G}$, $C_1 > 0$ and*

$$\lambda = \begin{cases} C_1 \Big( \dfrac{\log |\mathcal{G}|}{n} \Big)^{(q-1)/s}, & \text{when } q \leq s < 2q - 2, \\ C_1 \Big( \dfrac{\log |\mathcal{G}|}{n} \Big)^{q/(s+2)}, & \text{when } s \geq 2q - 2. \end{cases}$$

*The expected risk of the algorithm which draws uniformly its prediction function among $\mathbb{E}_{g \sim \pi_{-\lambda \Sigma_0}} g, \ldots, \mathbb{E}_{g \sim \pi_{-\lambda \Sigma_n}} g$ is upper bounded by*

$$\begin{cases} \min_{g \in \mathcal{G}} R(g) + C \Big( \dfrac{\log |\mathcal{G}|}{n} \Big)^{1 - (q-1)/s}, & \text{when } q \leq s \leq 2q - 2, \\ \min_{g \in \mathcal{G}} R(g) + C \Big( \dfrac{\log |\mathcal{G}|}{n} \Big)^{1 - q/s + 2}, & \text{when } s \geq 2q - 2, \end{cases}$$

*for a quantity $C$ which depends only on $C_1$, $b$, $A$, $q$ and $s$.*

PROOF. See Section 10.6. □

REMARK 7.1. In particular, for $q = 2$, with the minimal assumption $\mathbb{E}Y^2 \leq A$ (i.e., $s = 2$), the convergence rate is of order $n^{-1/2}$, and at the opposite, when $s$ goes to infinity, we recover the $n^{-1}$ rate we have under exponential moment condition on the output. Inequalities with precise constants for least square loss can also be found in the technical report [8], Section 7. For $q > 2$, low convergence rates (i.e., $n^{-\gamma}$ with $\gamma < 1/2$) appear when the moment assumption is weak: $\mathbb{E}|Y|^s \leq A$ for some $A > 0$ and $q \leq s < 2q - 2$. Convergence rates faster than the standard nonparametric rates $n^{-1/2}$ are achieved for $s > 2q - 2$. Fast convergence rates systematically occur when $1 < q < 2$ since for these values of $q$, we have $s \geq q > 2q - 2$. Surprisingly, for $q = 1$, the picture is completely different (see Section 8.3.2 for discussion and minimax optimality of the results of this section).

REMARK 7.2. Corollary 7.2 assumes that the prediction functions in $\mathcal{G}$ are uniformly bounded. It is an open problem to have the same kind of results under weaker assumptions such as a finite moment condition similar to the one used in Corollary 7.1.

**8. Lower bounds.** The simplest way to assess the quality of an algorithm and of its expected risk upper bound is to prove a risk lower bound saying that no algorithm has better convergence rate. This section provides this kind of assertion. The lower bounds developed here have the same spirit as the ones in [3, 14, 18], ([31], Chapter 15) and ([6], Section 5) to the extent that it relies on the following ideas:

- The supremum of a quantity $\mathcal{Q}(P)$ when the distribution $P$ belongs to some set $\mathcal{P}$ is larger than the supremum over a well-chosen finite subset of $\mathcal{P}$, and consequently is larger than the mean of $\mathcal{Q}(P)$ when the distribution $P$ is drawn uniformly in the finite subset.
- When the chosen subset is a hypercube of $2^m$ distributions (see Section 8.1), the design of a lower bound over the $2^m$ distributions reduces to the design of a lower bound over two distributions.
- When a data sequence $Z_1, \ldots, Z_n$ has similar likelihoods according to two different probability distributions, then no estimator will be accurate for both distributions: the maximum over the two distributions of the risk of any estimator trained on this sequence will be all the larger as the Bayes-optimal prediction associated with the two distributions are "far away."

We refer the reader to [15] and [46], Chapter 2, for lower bounds not particularly based on finding the appropriate hypercube. Our analysis focuses on hypercubes since in several settings they afford to obtain lower bounds with both the right convergence rate and close to optimal constants. Our contribution in this section is:

- to provide results for general nonregularized loss functions (we recall that nonregularized loss functions are loss functions which can be written as $L[(x,y), g] = \ell[y, g(x)]$ for some function $\ell: \mathcal{Y} \times \mathcal{Y} \to \mathbb{R}$),
- to improve the upper bound on the variational distance appearing in Assouad's argument,
- to generalize the argument to asymmetrical hypercubes which, to our knowledge, is the only way to find the lower bound matching the upper bound of Corollary 7.2 for $q \le s \le 2q-2$,
- to express the lower bounds in terms of similarity measures between two distributions characterizing the hypercube,
- to obtain lower bounds matching the upper bounds obtained in the previous sections.

REMARK 8.1. In [33], the optimality of the constant in front of the $(\log |\mathcal{G}|)/n$ has been proved by considering the situation when both $|\mathcal{G}|$ and $n$ go to infinity. Note that this worst-case analysis constant is not necessarily the same as our batch setting constant. This section shows that the batch setting constant is not "far" from the worst-case analysis constant.

Besides Lemma 4.3, which can be used to convert any worst-case analysis upper bounds into a risk upper bound in our batch setting, also means that any lower bounds for our batch setting lead to a lower bound in the sequential prediction setting (the converse is not true). Indeed the cumulative loss on the worst sequence of data is bigger than the average cumulative loss when the data are taken i.i.d. from some probability distribution. As a consequence, the bounds developed in this section partially solve the open problem introduced in [33], Section 3.4, consisting in developing tight nonasymptotical lower bounds. For least square loss and entropy loss, our bounds are off by a multiplicative factor smaller than 4 (see Remarks 8.5 and 8.4).

This section is organized as follows. Section 8.1 defines the quantities that characterize hypercubes of probability distributions and details the links between them. It also introduces a similarity measure between probability distributions coming from $f$-divergences (see [28]). We give our main lower bounds in Section 8.2. These bounds are illustrated in Section 8.3.

8.1. *Hypercube of probability distributions and f-similarities.*

DEFINITION 8.1. Let $m \in \mathbb{N}^*$. A hypercube of probability distributions is a family of $2^m$ probability distributions on $\mathcal{Z}$

$$\{P_{\bar{\sigma}} : \bar{\sigma} \triangleq (\sigma_1, \ldots, \sigma_m) \in \{-;+\}^m\}$$

having the same first marginal, denoted $\mu$,

$$P_{\bar{\sigma}}(dX) = P_{(+,\ldots,+)}(dX) \triangleq \mu(dX) \qquad \text{for any } \bar{\sigma} \in \{-;+\}^m,$$

and such that there exist:

- a partition $\mathcal{X}_0, \ldots, \mathcal{X}_m$ of $\mathcal{X}$ with $\mu(\mathcal{X}_1) = \cdots = \mu(\mathcal{X}_m)$,
- $h_1 \neq h_2$ in $\mathcal{Y}$,
- $0 \leq p_- < p_+ \leq 1$,

for which for any $j \in \{1, \ldots, m\}$, for any $x \in \mathcal{X}_j$, we have

$$(8.1) \qquad P_{\bar{\sigma}}(Y = h_1 | X = x) = p_{\sigma_j} = 1 - P_{\bar{\sigma}}(Y = h_2 | X = x),$$

and for any $x \in \mathcal{X}_0$, the distribution of $Y$ knowing $X = x$ is independent of $\bar{\sigma}$ (i.e., the $2^m$ conditional distributions are identical).

In particular, (8.1) means that for any $x \in \mathcal{X} - \mathcal{X}_0$, the conditional probability of the output knowing the input $x$ is concentrated on two values, and that, under the distribution $P_{\bar{\sigma}}$, the disproportion between the probabilities of these two values is all the larger as $p_{\sigma_j}$ is far from $1/2$ for $j$ the integer such that $x \in \mathcal{X}_j$.

An example of a hypercube is illustrated in Figure 3.

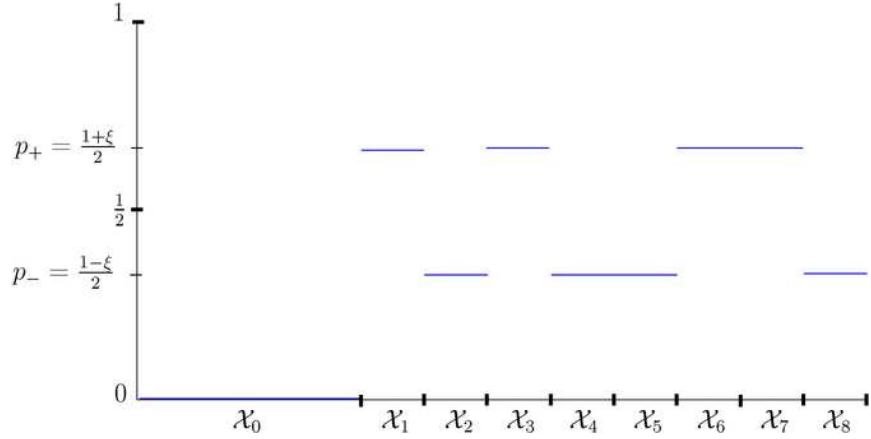

Fig. 3. *Representation of a probability distribution of the hypercube. Here the hypercube is symmetrical ($p_- = 1 - p_+$) with $m = 8$ and the probability distribution is characterized by $\bar{\sigma} = (+, -, +, -, -, +, +, -)$.*

REMARK 8.2. The use of hypercubes in which $p_+$ and $p_-$ are functions from $\mathcal{X} - \mathcal{X}_0$ to $[0; 1]$ and not just constants can be required when smoothness assumptions are put on the regression function $\eta : x \mapsto P(Y = 1 | X = x)$. This is typically the case in works on plug-in classifiers [2, 10]. For general hypercubes handling these kinds of constraints, we refer the reader to [8], Section 8.1.

Let $h_1$ and $h_2$ be distinct output values. For any $p \in [0; 1]$ and $y \in \mathcal{Y}$, consider

$$(8.2) \qquad \varphi_p(y) \triangleq p\ell(h_1, y) + (1 - p)\ell(h_2, y).$$

This is the risk of the prediction function identically equal to $y$ when the distribution generating the data satisfies $P[Y = y_1] = p = 1 - P[Y = y_2]$. Through this distribution, the quantity

$$(8.3) \qquad \phi(p) \triangleq \inf_{y \in \mathcal{Y}} \varphi_p(y)$$

can be viewed as the risk of the best constant prediction function.

For any $q_+$ and $q_-$ in $[0; 1]$, introduce

$$(8.4) \qquad \psi_{q_+, q_-}(\alpha) \triangleq \phi[\alpha q_+ + (1 - \alpha)q_-] - \alpha\phi(q_+) - (1 - \alpha)\phi(q_-).$$

DEFINITION 8.2. Let $\{P_{\bar{\sigma}} : \bar{\sigma} \triangleq (\sigma_1, \ldots, \sigma_m) \in \{-; +\}^m\}$ be a hypercube of distributions.

1. The positive integer $m$ is called the *dimension* of the hypercube.
2. The probability $w \triangleq \mu(\mathcal{X}_1) = \cdots = \mu(\mathcal{X}_m)$ is called the *edge probability.*
3. The *characteristic function of the hypercube* is the function $\tilde{\psi} : \mathbb{R}_+ \to \mathbb{R}_+$ defined as

$$(8.5) \qquad \tilde{\psi}(u) = \frac{mw}{2}(u+1)\psi_{p_+, p_-}\left(\frac{u}{u+1}\right).$$

4. The *edge discrepancy of type* I of the hypercube is

$$(8.6) \qquad d_{\mathrm{I}} \triangleq \frac{\tilde{\psi}(1)}{mw} = \psi_{p_+, p_-}(1/2)$$

5. The *edge discrepancy of type* II of the hypercube is defined as

$$(8.7) \qquad d_{\mathrm{II}} \triangleq \left(\sqrt{p_+(1-p_-)} - \sqrt{(1-p_+)p_-}\right)^2.$$

6. A probability distribution $P_0$ on $\mathcal{Z}$ satisfying $P_0(dX) = \mu(dX)$ and for any $x \in \mathcal{X} - \mathcal{X}_0$, $P_0[Y = h_1 | X = x] = \frac{1}{2} = P_0[Y = h_2 | X = x]$ will be referred to as a *base of the hypercube.*

7. Let $P_0$ be a base of the hypercube. Consider distributions $P_{[\sigma]}, \sigma \in \{-, +\}$ admitting the following density w.r.t. $P_0$:

$$\frac{P_{[\sigma]}}{P_0}(x, y) = \begin{cases} 2p_\sigma, & \text{when } x \in \mathcal{X}_1 \text{ and } y = h_1, \\ 2[1 - p_\sigma], & \text{when } x \in \mathcal{X}_1 \text{ and } y = h_2, \\ 1, & \text{otherwise.} \end{cases}$$

   The distributions $P_{[-]}$ and $P_{[+]}$ will be referred to as the *representatives of the hypercube.*

8. When the functions $p_+$ and $p_-$ satisfy $p_+ = 1 - p_-$ on $\mathcal{X} - \mathcal{X}_0$, the hypercube will be said to be *symmetrical.* In this case, the function $2p_+ - 1$ will be denoted $\xi$ so that

$$(8.8) \qquad \begin{aligned} p_+ &= \frac{1+\xi}{2}, \\[4pt] p_- &= \frac{1-\xi}{2}. \end{aligned}$$

   Otherwise it will be said to be *asymmetrical.*

9. A $(\tilde{m}, \tilde{w}, \tilde{d}_{\mathrm{II}})$-hypercube is a constant and symmetrical $\tilde{m}$-dimensional hypercube with edge probability $\tilde{w}$ and edge discrepancy of type II equal to $\tilde{d}_{\mathrm{II}}$.

Let us now give some properties of the quantities that have just been defined. The function $\phi$ is concave since it is the infimum of concave (affine) functions. Consequently, $\psi_{q_+, q_-}$ is concave and nonnegative on $[0; 1]$. Therefore $\tilde{\psi}$ is concave and nonnegative on $\mathbb{R}_+$ with $\tilde{\psi}(0) = 0$, hence $\tilde{\psi}$ is nondecreasing and satisfies

$$(8.9) \qquad \tilde{\psi}(u) \geq (u \wedge 1)\tilde{\psi}(1) = mw d_{\mathrm{I}}(u \wedge 1).$$

The edge discrepancies are both nonnegative quantities that are all the smaller as $p_-$ and $p_+$ become closer. When the function $\phi$ is twice differentiable on $]0;1[$, the edge discrepancy $d_{\mathrm{I}}$ can be written through

(8.10)
$$\psi_{p_+, p_-}(1/2)$$
$$= \frac{(p_+ - p_-)^2}{2} \int_0^1 [t \wedge (1-t)] |\phi''[tp_+ + (1-t)p_-]| \, dt,$$

which is proved by integration by parts.

For a $(\tilde{m}, \tilde{w}, \tilde{d_{\mathrm{II}}})$-hypercube, we have $m = \tilde{m}$, $w = \tilde{w}$, $d_{\mathrm{II}} = \tilde{d_{\mathrm{II}}}$, $\xi \equiv \sqrt{d_{\mathrm{II}}}$, $p_- \equiv (1 - \sqrt{d_{\mathrm{II}}})/2$ and $p_+ \equiv (1 + \sqrt{d_{\mathrm{II}}})/2$. So when $\phi$ is twice differentiable on $]p_-; p_+[$,

(8.11)
$$d_{\mathrm{I}} = \frac{d_{\mathrm{II}}}{2} \int_0^1 [t \wedge (1-t)] \left| \phi'' \left( \frac{1 - \sqrt{d_{\mathrm{II}}}}{2} + \sqrt{d_{\mathrm{II}}} t \right) \right| dt.$$

DEFINITION 8.3. When a probability distribution $\mathbb{P}$ is absolutely continuous w.r.t. another probability distribution $\mathbb{Q}$, that is, $\mathbb{P} \ll \mathbb{Q}$, $\frac{\mathbb{P}}{\mathbb{Q}}$ denotes the density of $\mathbb{P}$ w.r.t. $\mathbb{Q}$. Let $\mathbb{R}_+ = [0; +\infty[$. For any concave function $f : \mathbb{R}_+ \to \mathbb{R}_+$, we define the $f$-*similarity* between two probability distributions as

(8.12)
$$\mathcal{S}_f(\mathbb{P}, \mathbb{Q}) = \begin{cases} \int f\left(\dfrac{\mathbb{P}}{\mathbb{Q}}\right) d\mathbb{Q}, & \text{if } \mathbb{P} \ll \mathbb{Q}, \\ f(0), & \text{otherwise.} \end{cases}$$

We call it $f$-similarity in reference to $f$-divergence (see [28]) to which it is closely related. Here we use $f$-similarities since they are the quantities that naturally appear in our lower bounds.

8.2. *Generalized Assouad's lemma.* We recall that the $n$-fold product of a distribution $P$ is denoted $P^{\otimes n}$. We start this section with a general lower bound for hypercubes of distributions. This lower bound is expressed in terms of a similarity between $n$-fold products of representatives of the hypercube.

THEOREM 8.1. *Let $\mathcal{P}$ be a set of probability distributions containing a hypercube of distributions of characteristic function $\tilde{\psi}$ and representatives $P_{[-]}$ and $P_{[+]}$. For any training set size $n \in \mathbb{N}^*$ and any estimator $\hat{g}$, we have*

(8.13)
$$\sup_{P \in \mathcal{P}} \left\{ \mathbb{E} R(\hat{g}) - \min_g R(g) \right\} \geq \mathcal{S}_{\tilde{\psi}}(P_{[+]}^{\otimes n}, P_{[-]}^{\otimes n}),$$

*where the minimum is taken over the space of all prediction functions and $\mathbb{E} R(\hat{g})$ denotes the expected risk of the estimator $\hat{g}$ trained on a sample of size $n$: $\mathbb{E} R(\hat{g}) = \mathbb{E}_{Z_1^n \sim P^{\otimes n}} R(\hat{g}_{Z_1^n}) = \mathbb{E}_{Z_1^n \sim P^{\otimes n}} \mathbb{E}_{(X,Y) \sim P} \ell[Y, \hat{g}_{Z_1^n}(X)]$.*

Proof.   See Section 10.7.   □

This theorem provides a lower bound holding for any estimator and expressed in terms of the hypercube structure. To obtain a tight lower bound associated with a particular learning task, it then suffices to find the hypercube in $\mathcal{P}$ for which the r.h.s. of (8.13) is the largest possible. By providing lower bounds of $\mathcal{S}_{\bar{\psi}}(P_{[+]}^{\otimes n}, P_{[-]}^{\otimes n})$ that are more explicit w.r.t. the hypercube parameters, we obtain the following results that are more in a ready-to-use form than Theorem 8.1.

THEOREM 8.2.   *Let $\mathcal{P}$ be a set of probability distributions containing a hypercube of distributions characterized by its dimension $m$, its edge probability $w$ and its edge discrepancies $d_{\mathrm{I}}$ and $d_{\mathrm{II}}$ (see Definition 8.2). For any estimator $\hat{g}$ and training set size $n \in \mathbb{N}^*$, the following assertions hold:*

1. *We have*

$$\sup_{P \in \mathcal{P}} \left\{ \mathbb{E} R(\hat{g}) - \min_g R(g) \right\} \geq m w d_{\mathrm{I}} (1 - \sqrt{1 - [1 - d_{\mathrm{II}}]^{nw}})$$

(8.14)

$$\geq m w d_{\mathrm{I}} (1 - \sqrt{n w d_{\mathrm{II}}}).$$

2. *When the hypercube satisfies $p_+ \equiv 1 \equiv 1 - p_-$, we also have*

(8.15)          $$\sup_{P \in \mathcal{P}} \left\{ \mathbb{E} R(\hat{g}) - \min_g R(g) \right\} \geq m w d_{\mathrm{I}} (1 - w)^n.$$

Proof.   See Section 10.8.   □

The lower bound (8.15) is less general than (8.14) but provides results with tight constants when convergence rate of order $n^{-1}$ has to be proven (see Remarks 8.5 and 8.4).

REMARK 8.3.   The previous lower bounds consider deterministic estimators (or algorithms), that is, functions from the training set space $\bigcup_{n \geq 0} \mathcal{Z}^n$ to the prediction function space $\bar{\mathcal{G}}$. They still hold for randomized estimators, that is, functions from the training set space to the set $\mathcal{D}$ of probability distributions on $\bar{\mathcal{G}}$.

8.3. *Examples.*   Theorem 8.2 motivates the following simple strategy to obtain a lower bound for a given set $\mathcal{P}$ of probability distributions and a reference set $\mathcal{G}$ of prediction functions: it consists in looking for the hypercube contained in the set $\mathcal{P}$ and for which:

• the lower bound is maximized,

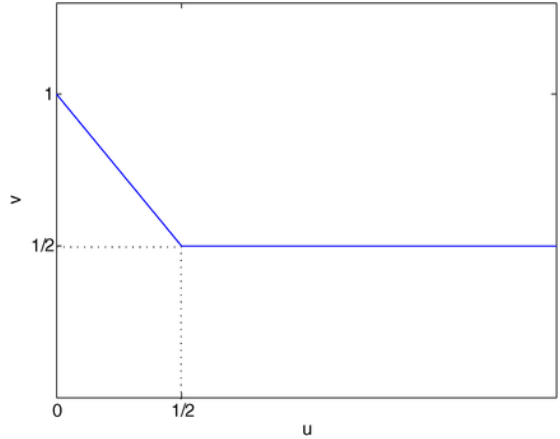

Fig. 4.   *Influence of the convexity of the loss on the optimal convergence rate. Let $c > 0$. We consider $L_q$-losses with $q = 1 + c(\frac{\log|\mathcal{G}|}{n})^u$ for $u \geq 0$. For such values of $q$, the optimal convergence rate of the associated learning task is of order $(\frac{\log|\mathcal{G}|}{n})^v$ with $1/2 \leq v \leq 1$. This figure represents the value of $u$ in abscissa and the value of $v$ in ordinate. The value $u = 0$ corresponds to constant $q$ greater than 1. For these $q$, the optimal convergence rate is of order $n^{-1}$ while for $q = 1$ or "very close" to 1, the convergence rate is of order $n^{-1/2}$.*

• for any distribution of the hypercube, $\mathcal{G}$ contains a best prediction function, that is, $\min_g R(g) = \min_{g \in \mathcal{G}} R(g)$.

In general, the order of the bound is given by the quantity $mwd_{\mathrm{I}}$ and the quantities $w$ and $d_{\mathrm{II}}$ are taken such that $nwd_{\mathrm{II}}$ is of order 1. This section illustrates this strategy by:

• providing learning lower bounds matching up to multiplicative constants the upper bounds developed in the previous sections,
• significantly improving the constants in classification lower bounds for Vapnik–Cervonenkis classes,
• showing that there is no uniform universal consistency for general loss functions.

8.3.1.  *$L_q$-regression with bounded outputs.*   We consider $\mathcal{Y} = [-B; B]$ and $\ell(y, y') = |y - y'|^q$, $q \geq 1$. The learning task is to predict as well as the best prediction function in a finite set $\mathcal{G}$ of cardinal denoted $|\mathcal{G}|$. The results of this section are roughly summed up in Figure 4, which represents the minimax optimal convergence rate for $L_q$-regression.

• *Case* $1 \leq q \leq 1 + \sqrt{\frac{\lfloor \log_2 |\mathcal{G}| \rfloor}{4n} \wedge 1}$. From (6.5), there exists an estimator $\hat{g}$ such that

(8.16)                    $\mathbb{E}R(\hat{g}) - \min_{g \in \mathcal{G}} R(g) \leq 2^{(2q-1)/2} B^q \sqrt{\frac{\log|\mathcal{G}|}{n}}$.

The following corollary of Theorem 8.2 shows that this result is tight.

**Theorem 8.3.** *Let $B > 0$ and $d \in \mathbb{N}^*$. For any training set size $n \in \mathbb{N}^*$ and any input space $\mathcal{X}$ containing at least $\lfloor \log_2 d \rfloor$ points, there exists a set $\mathcal{G}$ of $d$ prediction functions such that: for any estimator $\hat{g}$ there exists a probability distribution on the data space $\mathcal{X} \times [-B; B]$ for which*

$$\mathbb{E}R(\hat{g}) - \min_{g \in \mathcal{G}} R(g) \geq \begin{cases} c_q B^q \sqrt{\dfrac{\lfloor \log_2 |\mathcal{G}| \rfloor}{n}}, & \text{if } |\mathcal{G}| < 2^{4n+1}, \\ 2c_q B^q, & \text{otherwise,} \end{cases}$$

*where*

$$c_q = \begin{cases} \dfrac{1}{4}, & \text{if } q = 1, \\ \dfrac{q}{40}, & \text{if } 1 < q \leq 1 + \sqrt{\dfrac{\lfloor \log_2 |\mathcal{G}| \rfloor}{4n} \wedge 1}. \end{cases}$$

Proof. See Section 10.9.  □

*Case $q > 1 + \sqrt{\frac{\lfloor \log_2 |\mathcal{G}| \rfloor}{4n} \wedge 1}$.* We have seen in Section 4 that there exists an estimator $\hat{g}$ such that

$$(8.17) \qquad \mathbb{E}R(\hat{g}) - \min_{g \in \mathcal{G}} R(g) \leq \frac{q(1 \wedge 2^{q-2})B^q}{q-1}(\log 2)\frac{\log_2 |\mathcal{G}|}{n}.$$

The following corollary of Theorem 8.2 shows that this result is tight.

**Theorem 8.4.** *Let $B > 0$ and $d \in \mathbb{N}^*$. For any training set size $n \in \mathbb{N}^*$ and input space $\mathcal{X}$ containing at least $\lfloor \log_2(2d) \rfloor$ points, there exists a set $\mathcal{G}$ of $d$ prediction functions such that: for any estimator $\hat{g}$ there exists a probability distribution on the data space $\mathcal{X} \times [-B; B]$ for which*

$$\mathbb{E}R(\hat{g}) - \min_{g \in \mathcal{G}} R(g) \geq \left( \frac{q}{90(q-1)} \vee e^{-1} \right) B^q \left( \frac{\lfloor \log_2 |\mathcal{G}| \rfloor}{n+1} \wedge 1 \right).$$

Proof. See Section 10.9.  □

**Remark 8.4.** For least square regression (i.e., $q = 2$), Remark 8.5 holds provided that the multiplicative factor becomes $2e \log 2 \approx 3.77$. More generally, the method used here gives close to optimal constants but not the exact ones. We believe that this limit is due to the use of the hypercube structure. Indeed, the reader may check that for hypercubes of distributions, the upper bounds used in this section are not constant-optimal since the simplifying step consisting in using $\min_{\rho \in \mathcal{M}} \cdots \leq \min_{g \in \mathcal{G}} \cdots$ is loose.

The above analysis for $L_q$-losses can be generalized to show that there are essentialy two classes of bounded losses: the ones which are not convex or not enough convex (typical examples are the classification loss, the hinge loss and the absolute loss) and the ones which are sufficiently convex (typical examples are the least square loss, the entropy loss, the logit loss and the exponential loss). For the first class of losses, the edge discrepancy of type I is proportional to $\sqrt{d_{\mathrm{II}}}$ for constant and symmetrical hypercubes and (8.14) leads to a convergence rate of $\sqrt{(\log|\mathcal{G}|)/n}$. For the second class, the convergence rate is $(\log|\mathcal{G}|)/n$ and the lower bound can be explained by the fact that, when two prediction functions are different on a set with low probability (typically $n^{-1}$), it often happens that the training data have no input points in this set. For such training data, it is impossible to consistently choose the right prediction function.

This picture of convergence rates for finite models is rather well known, since:

- similar bounds (with looser constants) were known before for some cases (e.g., in classification; see [30, 50]),
- mutatis mutandis, the picture exactly matches the picture in the individual sequence prediction literature: for mixable loss functions (similar to "sufficiently convex"), the minimax regret is $O(\log|\mathcal{G}|)/n$, whereas for 0/1-type loss functions, it is $O(\sqrt{(\log|\mathcal{G}|)/n})$ (see, e.g., [33]).

8.3.2. $L_q$-regression for unbounded outputs having finite moments. The framework is similar to the one of Section 8.3.1 except that "$|Y| \leq B$ for some $B > 0$" is replaced with "$\mathbb{E}|Y|^s \leq A$ for some $s \geq q$ and $A > 0$."

*Case $q = 1$.* From (7.1), when $\sup_{g \in \mathcal{G}} E_Z g(X)^2 \leq b^2$ for some $b > 0$, there exists an estimator for which

$$\mathbb{E}R(\hat{g}) - \min_{g \in \mathcal{G}} R(g) \leq 2b\sqrt{(2\log|\mathcal{G}|)/n}.$$

The following corollary of Theorem 8.2 shows that this result is tight.

THEOREM 8.5. *For any training set size $n \in \mathbb{N}^*$, positive integer $d$, positive real number $b$ and input space $\mathcal{X}$ containing at least $\lfloor \log_2 d \rfloor$ points, there exists a set $\mathcal{G}$ of $d$ prediction functions uniformly bounded by $b$ such that: for any estimator $\hat{g}$ there exists a probability distribution for which $\mathbb{E}|Y| < +\infty$ and*

$$\mathbb{E}R(\hat{g}) - \min_{g \in \mathcal{G}} R(g) \geq \frac{b}{4}\sqrt{\frac{\lfloor \log_2|\mathcal{G}|\rfloor}{n} \wedge \frac{1}{4}}.$$

PROOF. Let $\tilde{m} = \lfloor \log_2|\mathcal{G}|\rfloor$. We consider a $(\tilde{m}, 1/\tilde{m}, \sqrt{\frac{\tilde{m}}{4n} \wedge 1})$-hypercube with $h_1 \equiv -b$ and $h_2 \equiv b$. One may check that $d_{\mathrm{I}} = b\sqrt{d_{\mathrm{II}}}$ so that (8.14) gives that for any estimator there exists a probability distribution for which $\mathbb{E}|Y| < +\infty$ and

$$\mathbb{E}R(\hat{g}) - \min_{g \in \mathcal{G}} R(g) \geq b\sqrt{\frac{\tilde{m}}{4n} \wedge 1}\left(1 - \sqrt{\frac{1}{4} \wedge \frac{n}{\tilde{m}}}\right),$$

hence the desired result.  □

*Case $q > 1$.* First let us recall the upper bound. In Corollary 7.2, under the assumptions

$$\begin{cases} \sup_{g \in \mathcal{G}, x \in \mathcal{X}} |g(x)| \leq b, & \text{for some } b > 0, \\ \mathbb{E}|Y|^s \leq A, & \text{for some } s \geq q \text{ and } A > 0, \\ \mathcal{G} \text{ finite}, \end{cases}$$

we have proposed an algorithm satisfying

$$R(\hat{g}) - \min_{g \in \mathcal{G}} R(g) \leq \begin{cases} C\left(\dfrac{\log|\mathcal{G}|}{n}\right)^{1-(q-1)/s}, & \text{when } q \leq s \leq 2q-2, \\ C\left(\dfrac{\log|\mathcal{G}|}{n}\right)^{1-q/(s+2)}, & \text{when } s \geq 2q-2, \end{cases}$$

for a quantity $C$ which depends only on $b$, $A$, $q$ and $s$.

The following corollary of Theorem 8.2 shows that this result is tight and is illustrated by Figure 5.

THEOREM 8.6.    *Let $d \in \mathbb{N}^*$, $s \geq q > 1$, $b > 0$ and $A > 0$. For any training set size $n \in \mathbb{N}^*$ and input space $\mathcal{X}$ containing at least $\lfloor \log_2(2d) \rfloor$ points, there exists a set $\mathcal{G}$ of $d$ prediction functions uniformly bounded by $b$ such that: for any estimator $\hat{g}$ there exists a probability distribution on the data space $\mathcal{X} \times \mathbb{R}$ for which $\mathbb{E}|Y|^s \leq A$ and*

$$\mathbb{E}R(\hat{g}) - \min_{g \in \mathcal{G}} R(g) \geq \begin{cases} C\left(\dfrac{\log|\mathcal{G}|}{n} \wedge 1\right)^{1-(q-1)/s}, \\ C\left(\dfrac{\log|\mathcal{G}|}{n} \wedge 1\right)^{1-q/(s+2)}, \end{cases}$$

*for a quantity $C$ which depends only on the real numbers $b$, $A$, $q$ and $s$.*

Both inequalities simultaneously hold but the first one is tight for $q \leq s \leq 2q - 2$ while the second one is tight for $s \geq 2q - 2$. They are both based on (8.14) applied to a $\lfloor \log_2|\mathcal{G}| \rfloor$-dimensional hypercube. Contrary to other lower bounds obtained in this work, the first inequality is based on asymmetrical hypercubes. The use of these kinds of hypercubes can be partially explained by the fact that the learning task is asymmetrical. Indeed all values of the output space do not have the same status since predictions are constrained to be in $[-b; b]$ while outputs are allowed to be in the whole real space (see the constraints on the hypercube in the proof given in Section 10.10).

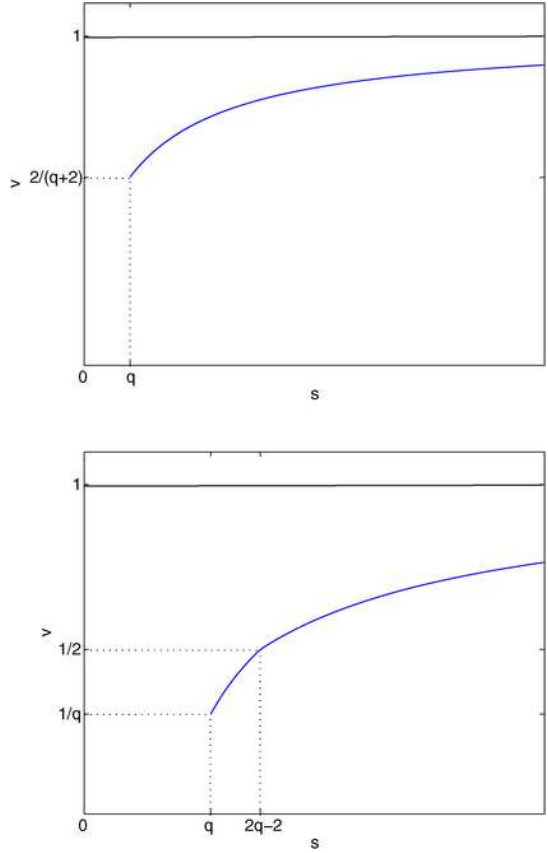

FIG. 5. *Optimal convergence rates in $L_q$-regression when the output has a finite moment of order $s$ (see Theorem 8.6). The convergence rate is of order $(\frac{\log|\mathcal{G}|}{n})^v$ with $0 < v \le 1$. The figure represents the value of $s$ in abscissa and the value of $v$ in ordinate. Two cases have to be distinguished. For $1 < q \le 2$ (figure on the top), $v$ depends smoothly on $q$. For $q > 2$ (figure on the bottom), two stages are observed depending whether $s$ is larger than $2q - 2$.*

8.3.3. *Entropy loss setting.* We consider $\mathcal{Y} = [0; 1]$ and $\ell(y, y') = K(y, y')$, where $K(y, y')$ is the Kullback–Leibler divergence between Bernoulli distributions with respective parameters $y$ and $y'$, that is, $K(y, y') = y \log(\frac{y}{y'}) + (1 - y) \log(\frac{1-y}{1-y'})$. We have seen in Section 4 that there exists an estimator $\hat{g}$ such that

$$(8.18) \qquad \mathbb{E}R(\hat{g}) - \min_{g \in \mathcal{G}} R(g) \le \frac{\log|\mathcal{G}|}{n}.$$

The following consequence of (8.15) shows that this result is tight.

THEOREM 8.7. *For any training set size $n \in \mathbb{N}^*$, positive integer $d$ and input space $\mathcal{X}$ containing at least $\lfloor \log_2(2d) \rfloor$ points, there exists a set $\mathcal{G}$ of $d$ prediction functions such that: for any estimator $\hat{g}$ there exists a probability distribution on the data space $\mathcal{X} \times [0; 1]$ for which*

$$\mathbb{E}R(\hat{g}) - \min_{g \in \mathcal{G}} R(g) \geq e^{-1}(\log 2)\left(1 \wedge \frac{\lfloor \log_2 |\mathcal{G}| \rfloor}{n+1}\right).$$

PROOF. We use a $(\tilde{m}, \frac{1}{n+1} \wedge \frac{1}{\tilde{m}}, 1)$-hypercube with $\tilde{m} = \lfloor \log_2 |\mathcal{G}| \rfloor = \lfloor \frac{\log |\mathcal{G}|}{\log 2} \rfloor$, $h_1 \equiv 0$ and $h_2 \equiv 1$. Let $H(y)$ denote the Shannon entropy of the Bernoulli distribution with parameter $y$, that is,

$$(8.19) \qquad H(y) = -y \log y - (1-y)\log(1-y).$$

Computations lead to: for any $p \in [0; 1]$,

$$\phi(p) = H(ph_1 + (1-p)h_2) - pH(h_1) - (1-p)H(h_2).$$

From (8.4) and Definition 8.2, we get

$$d_{\mathrm{I}} = \psi_{1,0,0,1}(1/2) = \phi_{0,1}(1/2) = H(1/2) = \log 2.$$

From (8.15), we obtain

$$\mathbb{E}R(\hat{g}) - \min_{g \in \mathcal{G}} R(g) \geq \left(\frac{\lfloor \log_2 |\mathcal{G}| \rfloor}{n+1} \wedge 1\right)(\log 2)\left(1 - \frac{1}{n+1} \wedge \frac{1}{\lfloor \log_2 |\mathcal{G}| \rfloor}\right)^n.$$

Then the result follows from $[1 - 1/(n+1)]^n \searrow e^{-1}$. □

REMARK 8.5. For $|\mathcal{G}| < 2^{n+2}$, the lower bound matches the upper bound (8.18) up to the multiplicative factor $e \approx 2.718$. For $|\mathcal{G}| \geq 2^{n+2}$, the size of the model is too large and, without any extra assumption, no estimator can learn from the data. To prove the result, we consider distributions for which the output is deterministic when knowing the input. So the lower bound does not come from noisy situations but from situations in which different prediction functions are not separated by the data to the extent that no input data fall into the (small) subset on which they are different.

8.3.4. *Binary classification.* We consider $\mathcal{Y} = \{0; 1\}$ and $l(y, y') = \mathbf{1}_{y \neq y'}$. Since the work of Vapnik and Cervonenkis [50], several lower bounds have been proposed and the most achieved ones are given in [30], Chapter 14. The following theorem provides an improvement of the constants of some of these bounds by a factor greater than 1000.

THEOREM 8.8. *Let $L \in [0; 1/2]$, $n \in \mathbb{N}$ and $\mathcal{G}$ be a set of prediction functions of VC-dimension $V \geq 2$. Consider the set $\mathcal{P}_L$ of probability distributions on $\mathcal{X} \times \{0; 1\}$ such that $\inf_{g \in \mathcal{G}} R(g) = L$. For any estimator $\hat{g}$:*

- *when $L = 0$, there exists $\mathbb{P} \in \mathcal{P}_0$ for which*

$$(8.20) \qquad \mathbb{E}R(\hat{g}) - \inf_{g \in \mathcal{G}} R(g) \geq \begin{cases} \dfrac{V-1}{2e(n+1)}, & \text{when } n \geq V-2, \\ \dfrac{1}{2}\left(1 - \dfrac{1}{V}\right)^n, & \text{otherwise,} \end{cases}$$

- *when $0 < L \leq 1/2$, there exists $\mathbb{P} \in \mathcal{P}_L$ for which*

$$\mathbb{E}R(\hat{g}) - \inf_{g \in \mathcal{G}} R(g)$$

$$(8.21)$$

$$\geq \begin{cases} \sqrt{\dfrac{L(V-1)}{32n}} \vee \dfrac{2(V-1)}{27n}, & \text{when } \dfrac{(1-2L)^2 n}{V} \geq \dfrac{4}{9}, \\ \dfrac{1-2L}{6}, & \text{otherwise,} \end{cases}$$

- *there exists a probability distribution for which*

$$(8.22) \qquad \mathbb{E}R(\hat{g}) - \inf_{g \in \mathcal{G}} R(g) \geq \dfrac{1}{8}\sqrt{\dfrac{V}{n}}.$$

SKETCH OF THE PROOF. For $h_1 \neq h_2$, we have $\phi(p) = p \wedge (1-p)$ and for symmetrical hypercubes $d_{\mathrm{I}} = \sqrt{d_{\mathrm{II}}}/2$. Then (8.20) comes from (8.15) and the use of a $(V-1, 1/(n+1), 1)$-hypercube and a $(V, 1/V, 1)$-hypercube.

To prove (8.21), from (8.14) and the use of a $(V-1, \frac{2L}{V-1}, \frac{V-1}{8nL})$-hypercube, a $(V-1, \frac{4}{9n(1-2L)^2}, (1-2L)^2)$-hypercube and a $(V, 1/V, (1-2L)^2)$-hypercube, we obtain

$$\mathbb{E}R(\hat{g}) - \inf_{g \in \mathcal{G}} R(g)$$

$$\geq \begin{cases} \sqrt{\dfrac{L(V-1)}{32n}}, & \text{when } \dfrac{(1-2L)^2 n}{V-1} \geq \dfrac{L}{2} \vee \dfrac{(1-2L)^2}{8L}, \\ \dfrac{2(V-1)}{27n(1-2L)}, & \text{when } \dfrac{(1-2L)^2 n}{V-1} \geq \dfrac{4}{9}, \\ \dfrac{1-2L}{2}\left(1 - \sqrt{\dfrac{(1-2L)^2 n}{V}}\right), & \text{always,} \end{cases}$$

which can be weakened into (8.21). Finally, (8.22) comes from the last inequality and by choosing $L$ such that $1 - 2L = \frac{1}{2}\sqrt{V/n}$. $\square$

In an asymptotical setting, [8], Section 8.4.3, provides a refinement of (8.22).

8.3.5. *No uniform universal consistency for general losses.* This type of result is well known and tells that there is no guarantee of doing well on finite samples. In a classification setting, when the input space is infinite, that is, $|\mathcal{X}| = +\infty$, by using a $(\lfloor n\alpha \rfloor, 1/\lfloor n\alpha \rfloor, 1)$-hypercube with $\alpha$ tending to infinity, one can recover that: for any training sample size $n$, "any discrimination rule can have an arbitrarily bad probability of error for finite sample size" [29], precisely:

$$\inf_{\hat{g}} \sup_{\mathbb{P}} \left\{ \mathbb{P}[Y \neq \hat{g}(X)] - \min_{g} \mathbb{P}[Y \neq g(X)] \right\} = 1/2,$$

where the infimum is taken over all (possibly randomized) classification rules. For general loss functions, as soon as $|\mathcal{X}| = +\infty$, we can use $(\lfloor n\alpha \rfloor, 1/\lfloor n\alpha \rfloor, 1)$-hypercubes with $\alpha$ tending to infinity and obtain

$$(8.23) \qquad \inf_{\hat{g}} \sup_{\mathbb{P}} \left\{ \mathbb{E}R(\hat{g}) - \inf_{g} R(g) \right\} \geq \sup_{y_1, y_2 \in \mathcal{Y}} \psi_{1,0,y_1,y_2}(1/2),$$

where $\psi$ is the function defined in (8.4).

**9. Summary of contributions and open problems.** This work has developed minimax optimal risk bounds for the general learning task consisting in predicting as well as the best function in a reference set. It has proposed to summarize this learning problem by the variance function appearing in the variance inequality (Section 3). The SeqRand algorithm (Figure 1) based on this variance function leads to minimax optimal convergence rates in the model selection aggregation problem, and our analysis gives a nice unified view to results coming from different communities.

In particular, results coming from the online learning literature are recovered in Section 4.1. The generalization error bounds obtained by Juditsky, Rigollet and Tsybakov in [34] are recovered for a slightly different algorithm in Section 5.

Without any extra assumption on the learning task, we have obtained a Bernstein's type bound which has no known equivalent form when the loss function is not assumed to be bounded (Section 6.1.1). When the loss function is bounded, the use of Hoeffding's inequality w.r.t. Gibbs distributions on the prediction function space instead of the distribution generating the data leads to an improvement by a factor 2 of the standard-style risk bound (Theorem 6.4).

To prove that our bounds are minimax optimal, we have refined Assouad's lemma particularly by taking into account the properties of the loss function. Theorem 8.2 is tighter than previous versions of Assouad's lemma and easier to apply to a learning setting than Fano's lemma (see, e.g., [46]); besides, the latter leads in general to very loose constants. It improves the constants of lower bounds related to Vapnik–Cervonenkis classes by a factor greater than 1000. We have also illustrated our upper and lower bounds by studying the influence of the noise of the output and of the convexity of the loss function.

For the $L_q$-loss with $q \geq 1$, new matching upper and lower bounds are given: in the online learning framework under boundedness assumption (Corollary 4.5 and Section 8.3.1 jointly with Remark 8.1), in the batch learning setting under boundedness assumption (Sections 4.1 and 8.3.1), in the batch learning setting for unbounded observations under moment assumptions (Sections 7 and 8.3.2). In the latter setting, we still do assume that the prediction functions are bounded. It is an open problem to replace this boundedness assumption with a moment condition.

Finally this work has the following limits. Most of our results concern expected risks and it is an open problem to provide corresponding tight exponential inequalities. Besides we should emphasize that our expected risk upper bounds hold only for our algorithm. This is quite different from the classical point of view that simultaneously gives upper bounds on the risk of any prediction function in the model. To our current knowledge, this classical approach has a flexibility that is not recovered in our approach. For instance, in several learning tasks, Dudley's chaining trick [32] is the only way to prove risk convergence with the optimal rate. So a natural question and another open problem is whether it is possible to combine the better variance control presented here with the chaining argument (or other localization argument used while exponential inequalities are available).

## 10. Proofs.

10.1. *Proof of Theorem 4.4.* First, by a scaling argument, it suffices to prove the result for $a = 0$ and $b = 1$. For $\mathcal{Y} = [0; 1]$, we modify the proof in Appendix A of [35]. Precisely, claims 1 and 2, with the notation used there, become:

1. If the function $f$ is concave in $\alpha([p; q])$, then we have $A_t(q) \leq B_t(p)$.
2. If $c \geq R(z, p, q)$ for any $z \in (p; q)$, then the function $f$ is concave in $\alpha([p; q])$.

Up to the missing $\alpha$ (typo), the difference is that we restrict ourselves to values of $z$ in $[p; q]$. The proof of claim 2 has no new argument. For claim 1, it suffices to modify the definition of $x'_{t,i}$ into $x'_{t,i} = q \wedge G^{-1}[\ell(p, x_{t,i})] \in [p; q]$. Then we have $L(p, x'_{t,i}) \leq L(p, x_{t,i})$ and $L(q, x'_{t,i}) \leq L(p, x_{t,i})$, hence $\alpha(x'_{t,i}) \geq \alpha(x_{t,i})$ and $\gamma(x'_{t,i}) \geq \gamma(x_{t,i})$. Now one can prove that $f$ is decreasing on $\alpha([p; q])$. By using Jensen's inequality, we get

$$\Delta_t(q) = -c \log \sum_{i=1}^{n} v_{t,i} \gamma(x_{t,i})$$

$$\geq -c \log \sum_{i=1}^{n} v_{t,i} \gamma(x'_{t,i})$$

$$= -c \log \sum_{i=1}^{n} v_{t,i} f[\alpha(x'_{t,i})]$$

$$\geq -c \log f\left[\sum_{i=1}^{n} v_{t,i} \alpha(x'_{t,i})\right]$$

$$\geq -c \log f\left[\sum_{i=1}^{n} v_{t,i} \alpha(x_{t,i})\right]$$

$$= L[q, G^{-1}(\Delta_t(p))].$$

The end of the proof of claim 1 is then identical.

10.2. *Proof of Theorem 6.1.* To check that the variance inequality holds, it suffices to prove that for any $z \in \mathcal{Z}$

$$(10.1) \qquad \mathbb{E}_{g' \sim \rho} \log \mathbb{E}_{g \sim \rho} e^{\lambda[L(z,g') - L(z,g)] - \lambda^2/2[L(z,g') - L(z,g)]^2} \leq 0.$$

To shorten formulae, let $\alpha(g', g) \triangleq \lambda[L(z, g') - L(z, g)]$. By Jensen's inequality and the following symmetrization trick, (10.1) holds:

$$\mathbb{E}_{g' \sim \rho} \mathbb{E}_{g \sim \rho} e^{\alpha(g',g) - \alpha^2(g',g)/2}$$

$$(10.2) \qquad \leq \tfrac{1}{2} \mathbb{E}_{g' \sim \rho} \mathbb{E}_{g \sim \rho} e^{\alpha(g',g) - \alpha^2(g',g)/2} + \tfrac{1}{2} \mathbb{E}_{g' \sim \rho} \mathbb{E}_{g \sim \rho} e^{-\alpha(g',g) - \alpha^2(g',g)/2}$$

$$\leq \mathbb{E}_{g' \sim \rho} \mathbb{E}_{g \sim \rho} \cosh(\alpha(g,g')) e^{-\alpha^2(g',g)/2} \leq 1,$$

where in the last inequality we used the inequality $\cosh(t) \leq e^{t^2/2}$ for any $t \in \mathbb{R}$. The result then follows from Theorem 3.1.

10.3. *Proof of Corollary 6.2.* To shorten the following formula, let $\mu$ denote the law of the prediction function produced by the SeqRand algorithm (w.r.t. simultaneously the training set and the randomizing procedure). Then (6.1) can be written as: for any $\rho \in \mathcal{M}$,

$$(10.3) \qquad \mathbb{E}_{g' \sim \mu} R(g') \leq \mathbb{E}_{g \sim \rho} R(g) + \frac{\lambda}{2} \mathbb{E}_{g \sim \rho} \mathbb{E}_{g' \sim \mu} V(g, g') + \frac{K(\rho, \pi)}{\lambda(n+1)}.$$

Define $\tilde{R}(g) = R(g) - R(\tilde{g})$ for any $g \in \mathcal{G}$. Under the generalized Mammen and Tsybakov assumption, for any $g, g' \in \mathcal{G}$, we have

$$\tfrac{1}{2} V(g, g') \leq \mathbb{E}_{Z \sim P}\{[L(Z, g) - L(Z, \tilde{g})]^2\} + \mathbb{E}_{Z \sim P}\{[L(Z, g') - L(Z, \tilde{g})]^2\}$$

$$\leq c\tilde{R}^{\gamma}(g) + c\tilde{R}^{\gamma}(g'),$$

so that (10.3) leads to

$$(10.4) \quad \mathbb{E}_{g' \sim \mu}[\tilde{R}(g') - c\lambda \tilde{R}^\gamma(g')] \leq \mathbb{E}_{g \sim \rho}[\tilde{R}(g) + c\lambda \tilde{R}^\gamma(g)] + \frac{K(\rho, \pi)}{\lambda(n+1)}.$$

This gives the first assertion. For the second statement, let $\tilde{u} \triangleq \mathbb{E}_{g' \sim \mu} \tilde{R}(g')$ and $\chi(u) \triangleq u - c\lambda u^\gamma$. By Jensen's inequality, the l.h.s. of (10.4) is lower bounded by $\chi(\tilde{u})$. By straightforward computations, for any $0 < \beta < 1$, when $u \geq (\frac{c\lambda}{1-\beta})^{1/(1-\gamma)}$, $\chi(u)$ is lower bounded by $\beta u$, which implies the desired result.

10.4. *Proof of Theorem 6.3.* Let us prove (6.3). Let $r(g)$ denote the empirical risk of $g \in \mathcal{G}$, that is, $r(g) = \frac{\Sigma_n(g)}{n}$. Let $\rho \in \mathcal{M}$ be some fixed distribution on $\mathcal{G}$. From [5], Section 8.1, with probability at least $1 - \varepsilon$ w.r.t. the training set distribution, for any $\mu \in \mathcal{M}$, we have

$$\mathbb{E}_{g' \sim \mu} R(g') - \mathbb{E}_{g \sim \rho} R(g)$$
$$\leq \mathbb{E}_{g' \sim \mu} r(g') - \mathbb{E}_{g \sim \rho} r(g) + \lambda \varphi(\lambda B) \mathbb{E}_{g' \sim \mu} \mathbb{E}_{g \sim \rho} V(g, g')$$
$$+ \frac{K(\mu, \pi) + \log(\varepsilon^{-1})}{\lambda n}.$$

Since the Gibbs distribution $\pi_{-\lambda \Sigma_n}$ minimizes $\mu \mapsto \mathbb{E}_{g' \sim \mu} r(g') + \frac{K(\mu, \pi)}{\lambda n}$, we have

$$\mathbb{E}_{g' \sim \pi_{-\lambda \Sigma_n}} R(g')$$
$$\leq \mathbb{E}_{g \sim \rho} R(g) + \lambda \varphi(\lambda B) \mathbb{E}_{g' \sim \pi_{-\lambda \Sigma_n}} \mathbb{E}_{g \sim \rho} V(g, g')$$
$$+ \frac{K(\rho, \pi) + \log(\varepsilon^{-1})}{\lambda n}.$$

Then we apply the following inequality:

$$\mathbb{E} W \leq \mathbb{E}(W \vee 0) = \int_0^{+\infty} \mathbb{P}(W > u) \, du = \int_0^1 \varepsilon^{-1} \mathbb{P}(W > \log(\varepsilon^{-1})) \, d\varepsilon$$

to the random variable

$$W = \lambda n [\mathbb{E}_{g' \sim \pi_{-\lambda \Sigma_n}} R(g') - \mathbb{E}_{g \sim \rho} R(g) - \lambda \varphi(\lambda B) \mathbb{E}_{g' \sim \pi_{-\lambda \Sigma_n}} \mathbb{E}_{g \sim \rho} V(g, g')]$$
$$- K(\rho, \pi).$$

We get $\mathbb{E} W \leq 1$. At last we may choose the distribution $\rho$ minimizing the upper bound to obtain (6.3). Similarly using [5], Section 8.3, we may prove (6.2).

10.5. *Proof of Lemma 6.5.* It suffices to apply the following adaptation of Lemma 5 of [55] to

$$\xi_i(Z_1, \ldots, Z_i) = L[Z_i, \mathcal{A}(Z_1^{i-1})] - L(Z_i, \tilde{g}).$$

LEMMA 10.1. *Let $\varphi$ still denote the positive convex increasing function defined as $\varphi(t) \triangleq \frac{e^t - 1 - t}{t^2}$. Let $b$ be a real number. For $i = 1, \ldots, n+1$, let $\xi_i \colon \mathcal{Z}^i \to \mathbb{R}$ be a function uniformly upper bounded by $b$. For any $\eta > 0$, $\varepsilon > 0$, with probability at least $1 - \varepsilon$ w.r.t. the distribution of $Z_1, \ldots, Z_{n+1}$, we have*

$$
\begin{aligned}
\sum_{i=1}^{n+1} \xi_i(Z_1, \ldots, Z_i) &\leq \sum_{i=1}^{n+1} \mathbb{E}_{Z_i} \xi_i(Z_1, \ldots, Z_i) \\
(10.5) \\
&\quad + \eta \varphi(\eta b) \sum_{i=1}^{n+1} \mathbb{E}_{Z_i} \xi_i^2(Z_1, \ldots, Z_i) + \frac{\log(\varepsilon^{-1})}{\eta},
\end{aligned}
$$

*where $\mathbb{E}_{Z_i}$ denotes the expectation w.r.t. the distribution of $Z_i$ only.*

REMARK 10.1. The same type of bounds without variance control can be found in [23].

PROOF OF LEMMA 10.1. For any $i \in \{0, \ldots, n+1\}$, define

$$\psi_i = \psi_i(Z_1, \ldots, Z_i) \triangleq \sum_{j=1}^{i} \xi_j - \sum_{j=1}^{i} \mathbb{E}_{Z_j} \xi_j - \eta \varphi(\eta b) \sum_{j=1}^{i} \mathbb{E}_{Z_j} \xi_j^2,$$

where $\xi_j$ is the short version of $\xi_j(Z_1, \ldots, Z_j)$. For any $i \in \{0, \ldots, n\}$, we trivially have

$$(10.6) \qquad \psi_{i+1} - \psi_i = \xi_{i+1} - \mathbb{E}_{Z_{i+1}} \xi_{i+1} - \eta \varphi(\eta b) \mathbb{E}_{Z_{i+1}} \xi_{i+1}^2.$$

Now for any $b \in \mathbb{R}$, $\eta > 0$ and any random variable $W$ such that $W \leq b$ a.s., we have

$$(10.7) \qquad\qquad \mathbb{E} e^{\eta(W - \mathbb{E}W - \eta \varphi(\eta b) \mathbb{E}W^2)} \leq 1.$$

REMARK 10.2. The proof of (10.7) is standard and can be found, for example, in [4], Section 7.1.1. We use (10.7) instead of the inequality used to prove Lemma 5 of [55], that is, $\mathbb{E} e^{\eta[W - \mathbb{E}W - \eta \varphi(\eta b') \mathbb{E}(W - \mathbb{E}W)^2]} \leq 1$ for $W - \mathbb{E}W \leq b'$ since we are interested in excess risk bounds. Precisely, we will take $W$ of the form $W = L(Z, g) - L(Z, g')$ for fixed functions $g$ and $g'$. Then we have $W \leq \sup_{z,g} L - \inf_{z,g} L$ while we only have $W - \mathbb{E}W \leq 2(\sup_{z,g} L - \inf_{z,g} L)$. Besides, the gain of having $\mathbb{E}(W - \mathbb{E}W)^2$ instead of $\mathbb{E}W^2$ is useless in the applications we develop here.

By combining (10.7) and (10.6), we obtain

$$\mathbb{E}_{Z_{i+1}} e^{\eta(\psi_{i+1} - \psi_i)} \leq 1. \tag{10.8}$$

By using Markov's inequality, we upper bound the following probability w.r.t. the distribution of $Z_1, \ldots, Z_{n+1}$:

$$\mathbb{P}\left( \sum_{i=1}^{n+1} \xi_i > \sum_{i=1}^{n+1} \mathbb{E}_{Z_i} \xi_i + \eta \varphi(\eta b) \sum_{i=1}^{n+1} \mathbb{E}_{Z_i} \xi_i^2 + \frac{\log(\varepsilon^{-1})}{\eta} \right)$$

$$= \mathbb{P}(\eta \psi_{n+1} > \log(\varepsilon^{-1}))$$

$$= \mathbb{P}(\varepsilon e^{\eta \psi_{n+1}} > 1)$$

$$\leq \varepsilon \mathbb{E} e^{\eta \psi_{n+1}}$$

$$\leq \varepsilon \mathbb{E}_{Z_1} \big( e^{\eta(\psi_1 - \psi_0)} \mathbb{E}_{Z_2} \big( \cdots e^{\eta(\psi_n - \psi_{n-1})} \mathbb{E}_{Z_{n+1}} e^{\eta(\psi_{n+1} - \psi_n)} \big) \big)$$

$$\leq \varepsilon,$$

where the last inequality follows from recursive use of (10.8).   □

### 10.6. *Proof of Corollary 7.2.*

We start with the following theorem concerning general loss functions.

**Theorem 10.2.** *Let $B \geq b > 0$ and $\mathcal{Y} = \mathbb{R}$. Consider a loss function $L$ which can be written as $L[(x,y),g] = \ell[y, g(x)]$, where the function $\ell : \mathbb{R} \times \mathbb{R} \to \mathbb{R}$ satisfies: there exists $\lambda_0 > 0$ such that for any $y \in [-B; B]$, the function $y' \mapsto e^{-\lambda_0 \ell(y,y')}$ is concave on $[-b; b]$. Let*

$$\Delta(y) = \sup_{|\alpha| \leq b, |\beta| \leq b} [\ell(y, \alpha) - \ell(y, \beta)].$$

*For $\lambda \in (0; \lambda_0]$, consider the algorithm that draws uniformly its prediction function in the set $\{\mathbb{E}_{g \sim \pi_{-\lambda \Sigma_0}} g, \ldots, \mathbb{E}_{g \sim \pi_{-\lambda \Sigma_n}} g\}$, and consider the deterministic version of this randomized algorithm. The expected risk of these algorithms satisfies*

$$\mathbb{E}_{Z_1^n} R\left( \frac{1}{n+1} \sum_{i=0}^{n} \mathbb{E}_{g \sim \pi_{-\lambda \Sigma_i}} g \right)$$

$$\leq \mathbb{E}_{Z_1^n} \frac{1}{n+1} \sum_{i=0}^{n} R(\mathbb{E}_{g \sim \pi_{-\lambda \Sigma_i}} g)$$

$$\leq \min_{\rho \in \mathcal{M}} \left\{ \mathbb{E}_{g \sim \rho} R(g) + \frac{K(\rho, \pi)}{\lambda(n+1)} \right\}$$

$$+ \mathbb{E} \left\{ \frac{\lambda \Delta^2(Y)}{2} \mathbf{1}_{\lambda \Delta(Y) < 1; |Y| > B} + \left[ \Delta(Y) - \frac{1}{2\lambda} \right] \mathbf{1}_{\lambda \Delta(Y) \geq 1; |Y| > B} \right\}. \tag{10.9}$$

Proof.  The first inequality follows from Jensen's inequality. Let us prove the second. According to Theorem 3.1, it suffices to check that the variance inequality holds for $0 < \lambda \le \lambda_0$, $\hat{\pi}(\rho)$ the Dirac distribution at $\mathbb{E}_{g \sim \rho} g$ and

$$\delta_\lambda[(x,y), g, g'] = \delta_\lambda(y) \triangleq \min_{0 \le \zeta \le 1} \left[ \zeta \Delta(y) + \frac{(1-\zeta)^2 \lambda \Delta^2(y)}{2} \right] \mathbf{1}_{|y|>B}$$

$$= \frac{\lambda \Delta^2(y)}{2} \mathbf{1}_{\lambda \Delta(y) < 1; |y|>B} + \left[ \Delta(y) - \frac{1}{2\lambda} \right] \mathbf{1}_{\lambda \Delta(y) \ge 1; |y|>B}.$$

- For any $z = (x,y) \in \mathcal{Z}$ such that $|y| \le B$, for any probability distribution $\rho$ and for the above values of $\lambda$ and $\delta_\lambda$, by Jensen's inequality, we have

$$\mathbb{E}_{g \sim \rho} e^{\lambda[L(z, \mathbb{E}_{g' \sim \rho} g') - L(z,g) - \delta_\lambda(z,g,g')]}$$

$$= e^{\lambda L(z, \mathbb{E}_{g' \sim \rho} g')} \mathbb{E}_{g \sim \rho} e^{-\lambda \ell[y, g(x)]}$$

$$\le e^{\lambda L(z, \mathbb{E}_{g' \sim \rho} g')} \left( \mathbb{E}_{g \sim \rho} e^{-\lambda_0 \ell[y, g(x)]} \right)^{\lambda/\lambda_0}$$

$$\le e^{\lambda \ell[y, \mathbb{E}_{g' \sim \rho} g'(x)] - \lambda \ell[y, \mathbb{E}_{g \sim \rho} g(x)]}$$

$$= 1,$$

where the last inequality comes from the concavity of $y' \mapsto e^{-\lambda_0 \ell(y, y')}$. This concavity argument goes back to [36], Section 4, and was also used in [19] and in some of the examples given in [34].

- For any $z = (x,y) \in \mathcal{Z}$ such that $|y| > B$, for any $0 \le \zeta \le 1$, by using twice Jensen's inequality and then by using the symmetrization trick presented in Section 6, we have

$$\mathbb{E}_{g \sim \rho} e^{\lambda[L(z, \mathbb{E}_{g' \sim \rho} g') - L(z,g) - \delta_\lambda(z,g,g')]}$$

$$= e^{-\delta_\lambda(y)} \mathbb{E}_{g \sim \rho} e^{\lambda[L(z, \mathbb{E}_{g' \sim \rho} g') - L(z,g)]}$$

$$\le e^{-\delta_\lambda(y)} \mathbb{E}_{g \sim \rho} e^{\lambda[\mathbb{E}_{g' \sim \rho} L(z,g') - L(z,g)]}$$

$$\le e^{-\delta_\lambda(y)} \mathbb{E}_{g \sim \rho} \mathbb{E}_{g' \sim \rho} e^{\lambda[L(z,g') - L(z,g)]}$$

$$= e^{-\delta_\lambda(y)} \mathbb{E}_{g \sim \rho} \mathbb{E}_{g' \sim \rho} \{ e^{\lambda(1-\zeta)[L(z,g') - L(z,g)] - 1/2 \lambda^2 (1-\zeta)^2 [L(z,g') - L(z,g)]^2}$$

$$\times e^{\lambda \zeta[L(z,g') - L(z,g)] + 1/2 \lambda^2 (1-\zeta)^2 [L(z,g') - L(z,g)]^2} \}$$

$$\le e^{-\delta_\lambda(y)} \mathbb{E}_{g \sim \rho} \mathbb{E}_{g' \sim \rho} \{ e^{\lambda(1-\zeta)[L(z,g') - L(z,g)] - 1/2 \lambda^2 (1-\zeta)^2 [L(z,g') - L(z,g)]^2}$$

$$\times e^{\lambda \zeta \Delta(y) + 1/2 \lambda^2 (1-\zeta)^2 \Delta^2(y)} \}$$

$$\le e^{-\delta_\lambda(y)} e^{\lambda \zeta \Delta(y) + (1/2) \lambda^2 (1-\zeta)^2 \Delta^2(y)}.$$

Taking $\zeta \in [0; 1]$ minimizing the last r.h.s., we obtain that

$$\mathbb{E}_{g \sim \rho} e^{\lambda[L(z, \mathbb{E}_{g' \sim \rho} g') - L(z,g) - \delta_\lambda(z,g,g')]} \le 1.$$

From the two previous computations, we obtain that for any $z \in \mathcal{Z}$,

$$\log \mathbb{E}_{g \sim \rho} e^{\lambda[L(z, \mathbb{E}_{g' \sim \rho} g') - L(z, g) - \delta_\lambda(z, g, g')]} \leq 0,$$

so that the variance inequality holds for the above values of $\lambda$, $\hat{\pi}(\rho)$ and $\delta_\lambda$, and the result follows from Theorem 3.1. $\square$

To apply Theorem 10.2, we will first determine $\lambda_0$ for which the function $\zeta : y' \mapsto e^{-\lambda_0|y-y'|^q}$ is concave. For any given $y \in [-B; B]$, for any $q > 1$, straightforward computations give

$$\zeta''(y') = [\lambda_0 q |y' - y|^q - (q-1)]\lambda_0 q |y' - y|^{q-2} e^{-\lambda_0|y-y'|^q}$$

for $y' \neq y$, hence $\zeta'' \leq 0$ on $[-b; b] - \{y\}$ for $\lambda_0 = \frac{q-1}{q(B+b)^q}$. Now since the derivative $\zeta'$ is defined at the point $y$, we conclude that the function $\zeta$ is concave on $[-b; b]$, so that we may use Theorem 10.2 with $\lambda_0 = \frac{q-1}{q(B+b)^q}$.

For any $|y| \geq b$, we have

$$2bq(|y| - b)^{q-1} \leq \Delta(y) \leq 2bq(|y| + b)^{q-1}.$$

As a consequence, when $|y| \geq b + (2bq\lambda)^{-1/(q-1)}$, we have $\lambda \Delta(y) \geq 1$ and $\Delta(y) - 1/(2\lambda)$ can be upper bounded by $C'|y|^{q-1}$, where the quantity $C'$ depends only on $b$ and $q$.

For other values of $|y|$, that is, when $b \leq |y| < b + (2bq\lambda)^{-1/(q-1)}$, we have

$$\begin{aligned}
&\frac{\lambda \Delta^2(y)}{2} \mathbf{1}_{\lambda \Delta(y) < 1; |y| > B} + \left[\Delta(y) - \frac{1}{2\lambda}\right] \mathbf{1}_{\lambda \Delta(y) \geq 1; |y| > B} \\
&= \min_{0 \leq \zeta \leq 1} \left[\zeta \Delta(y) + \frac{(1 - \zeta)^2 \lambda \Delta^2(y)}{2}\right] \mathbf{1}_{|y| > B} \\
&\leq \frac{1}{2} \lambda \Delta^2(y) \mathbf{1}_{|y| > B} \\
&\leq 2\lambda b^2 q^2 (|y| + b)^{2q-2} \mathbf{1}_{|y| > B} \\
&\leq C'' \lambda |y|^{2q-2} \mathbf{1}_{|y| > B},
\end{aligned}$$

where $C''$ depends only on $b$ and $q$.

Therefore, from (10.9), for any $0 < b \leq B$ and $\lambda > 0$ satisfying $\lambda \leq \frac{q-1}{q(B+b)^q}$, the expected risk is upper bounded by

$$\begin{aligned}
(10.10) \quad &\min_{\rho \in \mathcal{M}} \left\{ \mathbb{E}_{g \sim \rho} R(g) + \frac{K(\rho, \pi)}{\lambda(n+1)} \right\} + \mathbb{E}\left\{ C'|Y|^{q-1} \mathbf{1}_{|Y| \geq b + (2bq\lambda)^{-1/(q-1)}; |Y| > B} \right\} \\
&+ \mathbb{E}\left\{ C'' \lambda |Y|^{2q-2} \mathbf{1}_{B < |Y| < b + (2bq\lambda)^{-1/(q-1)}} \right\}.
\end{aligned}$$

Let us take $B = (\frac{q-1}{q\lambda})^{1/q} - b$ with $\lambda$ small enough to ensure that $b \leq B \leq b + (2bq\lambda)^{-1/(q-1)}$. This means that $\lambda$ should be taken smaller than some positive constant depending only on $b$ and $q$. Then (10.10) can be written as

$$\min_{\rho \in \mathcal{M}} \left\{ \mathbb{E}_{g \sim \rho} R(g) + \frac{K(\rho, \pi)}{\lambda(n+1)} \right\} + \mathbb{E}\{C'|Y|^{q-1} \mathbf{1}_{|Y| \geq b + (2bq\lambda)^{-1/(q-1)}}\}$$

$$+ \mathbb{E}\{C''\lambda |Y|^{2q-2} \mathbf{1}_{((q-1)/(q\lambda))^{1/q} - b < |Y| < b + (2bq\lambda)^{-1/(q-1)}}\}.$$

The moment assumption on $Y$ implies

(10.11)     $\alpha^{s-q} \mathbb{E}|Y|^q \mathbf{1}_{|Y| \geq \alpha} \leq A$     for any $0 \leq q \leq s$ and $\alpha \geq 0$.

So we can upper bound (10.10) with

$$\min_{g \in \mathcal{G}} R(g) + \frac{\log |\mathcal{G}|}{\lambda n} + C\lambda^{(s+1-q)/(q-1)}$$

$$+ C\lambda(\lambda^{(s-2q+2)/q} \mathbf{1}_{s \geq 2q-2} + \lambda^{(2-2q+s)(q-1)} \mathbf{1}_{s<2q-2}),$$

where $C$ depends only on $b$, $A$, $q$ and $s$. So we get

$$\mathbb{E}_{Z_1^n} \frac{1}{n+1} \sum_{i=0}^{n} R(\mathbb{E}_{g \sim \pi_{-\lambda \Sigma_i}} g)$$

$$\leq \min_{g \in \mathcal{G}} R(g) + \frac{\log |\mathcal{G}|}{\lambda n} + C\lambda^{(s+1-q)/(q-1)} + C\lambda^{(s-q+2)/q} \mathbf{1}_{s \geq 2q-2}$$

$$\leq \min_{g \in \mathcal{G}} R(g) + \frac{\log |\mathcal{G}|}{\lambda n} + C\lambda^{(s+1-q)/(q-1)} \mathbf{1}_{s<2q-2} + C\lambda^{(s-q+2)/q} \mathbf{1}_{s \geq 2q-2},$$

since $\frac{s+1-q}{q-1} \geq \frac{s-q+2}{q}$ is equivalent to $s \geq 2q-2$. By taking $\lambda$ of order of the minimum of the r.h.s. (which implies that $\lambda$ goes to 0 when $n/\log |\mathcal{G}|$ goes to infinity), we obtain the desired result.

10.7. *Proof of Theorem 8.1.* The symbols $\sigma_1, \ldots, \sigma_m$ still denote the coordinates of $\bar{\sigma} \in \{-; +\}^m$. For any $r \in \{-; 0; +\}$, define

$$\bar{\sigma}_{j,r} \triangleq (\sigma_1, \ldots, \sigma_{j-1}, r, \sigma_{j+1}, \ldots, \sigma_m)$$

as the vector deduced from $\bar{\sigma}$ by fixing its $j$th coordinate to $r$. Since $\bar{\sigma}_{j,+}$ and $\bar{\sigma}_{j,-}$ belong to $\{-; +\}^m$, we have already defined $P_{\bar{\sigma}_{j,+}}$ and $P_{\bar{\sigma}_{j,-}}$. Now we define the distribution $P_{\bar{\sigma}_{j,0}}$ as $P_{\bar{\sigma}_{j,0}}(dX) = \mu(dX)$ and

$$1 - P_{\bar{\sigma}_{j,0}}(Y=h_2|X) = P_{\bar{\sigma}_{j,0}}(Y=h_1|X)$$

$$= \begin{cases} \frac{1}{2}, & \text{for any } X \in \mathcal{X}_j, \\ P_{\bar{\sigma}}(Y=h_1|X), & \text{otherwise.} \end{cases}$$

The distribution $P_{\bar{\sigma}_{j,0}}$ differs from $P_{\bar{\sigma}}$ only by the conditional law of the output knowing that the input is in $\mathcal{X}_j$. We recall that $P^{\otimes n}$ denotes the $n$-fold product of a distribution $P$. For any $r \in \{-;+\}$, introduce the likelihood ratios for the data $Z_1^n = (Z_1, \ldots, Z_n)$: $\pi_{r,j}(Z_1^n) \triangleq \frac{P_{\bar{\sigma}_{j,r}}^{\otimes n}}{P_{\bar{\sigma}_{j,0}}^{\otimes n}}(Z_1^n)$. This quantity is independent of the value of $\bar{\sigma}$. Let $\nu$ be the uniform distribution on $\{-,+\}$, that is, $\nu(\{+\}) = 1/2 = 1 - \nu(\{-\})$. In the following, $\mathbb{E}_{\bar{\sigma}}$ denotes the expectation when $\bar{\sigma}$ is drawn according to the $m$-fold product distribution of $\nu$, and $\mathbb{E}_X = \mathbb{E}_{X \sim \mu}$. We have

$$\sup_{P \in \mathcal{P}} \left\{ \mathbb{E}_{Z_1^n \sim P^{\otimes n}} R(\hat{g}) - \min_g R(g) \right\}$$

$$\geq \sup_{\bar{\sigma} \in \{-;+\}^m} \left\{ \mathbb{E}_{Z_1^n \sim P_{\bar{\sigma}}^{\otimes n}} \mathbb{E}_{Z \sim P_{\bar{\sigma}}} \ell[Y, \hat{g}(X)] - \min_g \mathbb{E}_{Z \sim P_{\bar{\sigma}}} \ell[Y, g(X)] \right\}$$

$$= \sup_{\bar{\sigma} \in \{-;+\}^m} \left\{ \mathbb{E}_{Z_1^n \sim P_{\bar{\sigma}}^{\otimes n}} \mathbb{E}_{X \sim P_{\bar{\sigma}}(dX)} \left[ \mathbb{E}_{Y \sim P_{\bar{\sigma}}(dY|X)} \ell[Y, \hat{g}(X)] \right. \right.$$
$$\left. \left. - \min_{y \in \mathcal{Y}} \mathbb{E}_{Y \sim P_{\bar{\sigma}}(dY|X)} \ell(Y, y) \right] \right\}$$

$$(10.12) \quad = \sup_{\bar{\sigma} \in \{-;+\}^m} \left\{ \mathbb{E}_{Z_1^n \sim P_{\bar{\sigma}}^{\otimes n}} \mathbb{E}_X \left[ \sum_{j=0}^m \mathbf{1}_{X \in \mathcal{X}_j} (\varphi_{p_{\sigma_j}}[\hat{g}(X)] - \phi[p_{\sigma_j}]) \right] \right\}$$

$$\geq \mathbb{E}_{\bar{\sigma}} \mathbb{E}_{Z_1^n \sim P_{\bar{\sigma}}^{\otimes n}} \mathbb{E}_X \left[ \sum_{j=1}^m \mathbf{1}_{X \in \mathcal{X}_j} (\varphi_{p_{\sigma_j}}[\hat{g}(X)] - \phi[p_{\sigma_j}]) \right]$$

$$= \sum_{j=1}^m \mathbb{E}_X \left\{ \mathbf{1}_{X \in \mathcal{X}_j} \mathbb{E}_{\bar{\sigma}} \mathbb{E}_{Z_1^n \sim P_{\bar{\sigma}_{j,0}}^{\otimes n}} \left[ \frac{P_{\bar{\sigma}}^{\otimes n}}{P_{\bar{\sigma}_{j,0}}^{\otimes n}} (Z_1^n) (\varphi_{p_{\sigma_j}}[\hat{g}(X)] - \phi[p_{\sigma_j}]) \right] \right\}$$

$$= \sum_{j=1}^m \mathbb{E}_X \left\{ \mathbf{1}_{X \in \mathcal{X}_j} \mathbb{E}_{\sigma_1, \ldots, \sigma_{j-1}, \sigma_{j+1}, \ldots, \sigma_m} \mathbb{E}_{Z_1^n \sim P_{\bar{\sigma}_{j,0}}^{\otimes n}} \mathbb{E}_{\sigma_j \sim \nu} \pi_{\sigma_j, j}(Z_1^n) \right.$$
$$\left. \times (\varphi_{p_{\sigma_j}}[\hat{g}(X)] - \phi[p_{\sigma_j}]) \right\}.$$

The two inequalities in (10.12) are Assouad's argument [3]. For any $x \in \mathcal{X}$, introduce $\alpha_j(Z_1^n) = \frac{\pi_{+,j}(Z_1^n)}{\pi_{+,j}(Z_1^n) + \pi_{-,j}(Z_1^n)}$. The last expectation in (10.12) is

$$\mathbb{E}_{\sigma \sim \nu} \pi_{\sigma,j}(Z_1^n) (\varphi_{p_\sigma}[\hat{g}(X)] - \phi[p_\sigma])$$

$$= \frac{1}{2} [\pi_{+,j}(Z_1^n) + \pi_{-,j}(Z_1^n)]$$
$$\times \{ \alpha_j(Z_1^n) \varphi_{p_+}[\hat{g}(X)] + [1 - \alpha_j(Z_1^n)] \varphi_{p_-}[\hat{g}(X)]$$
$$(10.13) \qquad - \alpha_j(Z_1^n) \phi(p_+) - [1 - \alpha_j(Z_1^n)] \phi(p_-) \}$$

$$= \frac{1}{2} [\pi_{+,j}(Z_1^n) + \pi_{-,j}(Z_1^n)] \{ \varphi_{\alpha_j(Z_1^n)p_+ + [1 - \alpha_j(Z_1^n)]p_-} [\hat{g}(X)]$$

$$-\alpha_j(Z_1^n)\phi(p_+)-[1-\alpha_j(Z_1^n)]\phi(p_-)\}$$

$$\geq \frac{1}{2}[\pi_{+,j}(Z_1^n)+\pi_{-,j}(Z_1^n)]\{\phi(\alpha_j(Z_1^n)p_+ + [1-\alpha_j(Z_1^n)]p_-$$
$$-\alpha_j(Z_1^n)\phi(p_+)-[1-\alpha_j(Z_1^n)]\phi(p_-)\}$$

$$=\frac{1}{2}[\pi_{+,j}(Z_1^n)+\pi_{-,j}(Z_1^n)]\psi_{p_+,p_-}[\alpha_j(Z_1^n)]$$

$$=\frac{1}{mw}\pi_{-,j}(Z_1^n)\tilde{\psi}\left(\frac{\pi_{+,j}(Z_1^n)}{\pi_{-,j}(Z_1^n)}\right)$$

so that

$$\sup_{P\in\mathcal{P}}\left\{\mathop{\mathbb{E}}_{Z_1^n\sim P^{\otimes n}}R(\hat{g})-\min_g R(g)\right\}$$

$$\geq \frac{1}{mw}\sum_{j=1}^m \mathbb{E}_X\left\{\mathbf{1}_{X\in\mathcal{X}_j}\mathbb{E}_{\bar{\sigma}}\mathbb{E}_{Z_1^n\sim P^{\otimes n}_{\bar{\sigma},0}}\left[\pi_{-,j}(Z_1^n)\tilde{\psi}\left(\frac{\pi_{+,j}(Z_1^n)}{\pi_{-,j}(Z_1^n)}\right)\right]\right\}$$

$$=\frac{1}{mw}\sum_{j=1}^m \mathbb{E}_X\left\{\mathbf{1}_{X\in\mathcal{X}_j}\mathbb{E}_{\bar{\sigma}}\mathcal{S}_{\tilde{\psi}}(P^{\otimes n}_{\bar{\sigma}_{j,+}},P^{\otimes n}_{\bar{\sigma}_{j,-}})\right\}$$

$$=\frac{1}{m}\sum_{j=1}^m \mathbb{E}_{\bar{\sigma}}\mathcal{S}_{\tilde{\psi}}(P^{\otimes n}_{\bar{\sigma}_{j,+}},P^{\otimes n}_{\bar{\sigma}_{j,-}}).$$

Now since we consider a hypercube, for any $j\in\{1,\dots,m\}$, all the terms in the sum are equal. Besides one can check that the last $f$-similarity does not depend on $\bar{\sigma}$, and is equal to $\mathcal{S}_{\tilde{\psi}}(P^{\otimes n}_{[+]},P^{\otimes n}_{[-]})$ where we recall that $P_{[+]}$ and $P_{[-]}$ denote the representatives of the hypercube (see Definition 8.2) Therefore we obtain

$$\sup_{P\in\mathcal{P}}\left\{\mathbb{E}R(\hat{g})-\min_g R(g)\right\}\geq \mathcal{S}_{\tilde{\psi}}(P^{\otimes n}_{[+]},P^{\otimes n}_{[-]}).$$

10.8. *Proof of Theorem 8.2.* First, when the hypercube satisfies $p_+ = 1 = 1-p_-$, from the definition of $d_{\mathrm{I}}$ given in (8.6), we have $\mathcal{S}_{\tilde{\psi}}(P^{\otimes n}_{[+]},P^{\otimes n}_{[-]}) = mwd_{\mathrm{I}}(1-w)^n$ so that Theorem 8.1 implies (8.15).

Inequality (8.14) is deduced from Theorem 8.1 by lower bounding the $\tilde{\psi}$-similarity. Since $u\mapsto u\wedge 1$ is a nonnegative concave function defined on $\mathbb{R}_+$, we may define the similarity

$$\mathcal{S}_{\wedge}(\mathbb{P},\mathbb{Q})\triangleq \int\left(\frac{\mathbb{P}}{\mathbb{Q}}\wedge 1\right)d\mathbb{Q}=\int(d\mathbb{P}\wedge d\mathbb{Q}),$$

where the second equality introduces a formal (but intuitive) notation. From Theorem 8.1, by using (8.9), we obtain

COROLLARY 10.3.   *Let $\mathcal{P}$ be a set of probability distributions containing a hypercube of distributions of characteristic function $\tilde{\psi}$ and representatives $P_{[-]}$ and $P_{[+]}$. For any estimator $\hat{g}$, we have*

$$(10.14) \qquad \sup_{P \in \mathcal{P}} \left\{ \mathbb{E} R(\hat{g}) - \min_g R(g) \right\} \geq m w d_{\mathrm{I}} \mathcal{S}_\wedge (P_{[+]}^{\otimes n}, P_{[-]}^{\otimes n}),$$

*where the minimum is taken over the space of prediction functions.*

The following lemma and (10.14) imply (8.14).

LEMMA 10.4.   *We have*

$$(10.15) \qquad \mathcal{S}_\wedge (P_{[+]}^{\otimes n}, P_{[-]}^{\otimes n}) \geq 1 - \sqrt{1 - [1 - d_{\mathrm{II}}]^{nw}} \geq 1 - \sqrt{n w d_{\mathrm{II}}}.$$

PROOF.   See Section 10.8.1.   □

10.8.1. *Proof of Lemma 10.4.*   For $\sigma \in \{-, +\}$, define $Q_\sigma$ as the probability on $\{h_1, h_2\}$ such that $Q_\sigma(Y = h_1) = p_\sigma = 1 - Q_\sigma(Y = h_2)$. The following lemma relates the $\wedge$-similarity between representatives of the hypercube and the $\wedge$-similarity between $Q_+$ and $Q_-$.

LEMMA 10.5.   *Consider a convex function $\gamma : \mathbb{R}_+ \to \mathbb{R}_+$ such that*

$$\gamma(k) \leq \mathcal{S}_\wedge (Q_+^{\otimes k}, Q_-^{\otimes k})$$

*for any $k \in \{0, \ldots, n\}$, where by convention $\mathcal{S}_\wedge(Q_+^{\otimes 0}, Q_-^{\otimes 0}) = 1$. For any estimator $\hat{g}$, we have*

$$\mathcal{S}_\wedge (P_{[+]}^{\otimes n}, P_{[-]}^{\otimes n}) \geq \gamma(nw).$$

PROOF.   For any points $z_1 = (x_1, y_1), \ldots, z_n = (x_n, y_n)$ in $\mathcal{X} \times \{h_1, h_2\}$, let $\mathcal{C}(z_1, \ldots, z_n)$ denote the number of $z_i$ for which $x_i \in \mathcal{X}_1$. For any $k \in \{0, \ldots, n\}$, let $B_k = \mathcal{C}^{-1}(\{k\})$ denote the subset of $(\mathcal{X} \times \{h_1, h_2\})^n$ for which exactly $k$ points are in $\mathcal{X}_1 \times \{h_1, h_2\}$. We recall that there are $\binom{n}{k}$ possibilities of taking $k$ elements among $n$ and the probability of $X \in \mathcal{X}_1$ when $X$ is drawn according to $\mu$ is $w = \mu(\mathcal{X}_1)$. Let $\mathcal{Z}_1 = \mathcal{X}_1 \times \{h_1, h_2\}$ and let $\mathcal{Z}_1^c$ denote the complement of $\mathcal{Z}_1$. We have

$$\mathcal{S}_\wedge (P_{[+]}^{\otimes n}, P_{[-]}^{\otimes n})$$

$$= \int 1 \wedge \left( \frac{P_{[+]}^{\otimes n}}{P_{[-]}^{\otimes n}}(z_1, \ldots, z_n) \right) dP_{[-]}^{\otimes n}(z_1, \ldots, z_n)$$

$$(10.16) \quad = \sum_{k=0}^{n} \int_{B_k} 1 \wedge \left( \frac{P_{[+]}}{P_{[-]}}(z_1) \cdots \frac{P_{[+]}}{P_{[-]}}(z_n) \right) dP_{[-]}(z_1) \cdots dP_{[-]}(z_n)$$

$$= \sum_{k=0}^{n} \binom{n}{k} \int_{(\mathcal{Z}_1)^k \times (\mathcal{Z}_1^c)^{n-k}} 1 \wedge \left( \frac{P_{[+]}}{P_{[-]}}(z_1) \cdots \frac{P_{[+]}}{P_{[-]}}(z_n) \right) dP_{[-]}^{\otimes n}(z_1, \ldots, z_n)$$

$$= \sum_{k=0}^{n} \binom{n}{k} \int_{(\mathcal{Z}_1)^k \times (\mathcal{Z}_1^c)^{n-k}} 1 \wedge \left( \frac{P_{[+]}}{P_{[-]}}(z_1) \cdots \frac{P_{[+]}}{P_{[-]}}(z_k) \right) dP_{[-]}^{\otimes n}(z_1, \ldots, z_n)$$

$$= \sum_{k=0}^{n} \binom{n}{k} \mu^{n-k}(\mathcal{Z}_1^c)$$

$$\times \int_{(\mathcal{Z}_1)^k} 1 \wedge \left( \frac{P_{[+]}}{P_{[-]}}(z_1) \cdots \frac{P_{[+]}}{P_{[-]}}(z_k) \right) dP_{[-]}^{\otimes n}(z_1, \ldots, z_k)$$

$$= \sum_{k=0}^{n} \binom{n}{k} \mu^{n-k}(\mathcal{Z}_1^c) \mu^{k}(\mathcal{Z}_1) \mathcal{S}_{\wedge}(Q_+^{\otimes k}, Q_-^{\otimes k})$$

$$\geq \sum_{k=0}^{n} \binom{n}{k} (1-w)^{n-k} w^k \gamma(k)$$

$$= \mathbb{E} \gamma(V),$$

where $V$ is a Binomial distribution with parameters $n$ and $w$. By Jensen's inequality, we have $\mathbb{E}\gamma(V) \geq \gamma[\mathbb{E}(V)] = \gamma(nw)$, which ends the proof. $\square$

The interest of the previous lemma is to provide a lower bound on the similarity between representatives of the hypercube from a lower bound on the similarities between distributions much simpler to study. The following result lower bounds these similarities.

LEMMA 10.6.  *For any nonnegative integer $k$, we have*

$$(10.17) \qquad \mathcal{S}_{\wedge}(Q_+^{\otimes k}, Q_-^{\otimes k}) \geq 1 - \sqrt{1 - [1 - d_{\mathrm{II}}]^k} \geq 1 - \sqrt{k d_{\mathrm{II}}}.$$

PROOF.  To study divergences (or equivalently, similarities) between $k$-fold product distributions, the standard way is to link the divergence (or similarity) of the product with the ones of base distributions. This leads to tensorization equalities or inequalities. To obtain a tensorization inequality for $\mathcal{S}_{\wedge}$, we introduce the similarity associated with the square root function (which is nonnegative and concave):

$$\mathcal{S}_{\sqrt{}}(\mathbb{P}, \mathbb{Q}) \triangleq \int \sqrt{d\mathbb{P}\, d\mathbb{Q}}$$

and use the following lemmas:

LEMMA 10.7.  *For any probability distributions $\mathbb{P}$ and $\mathbb{Q}$, we have*

$$\mathcal{S}_{\wedge}(\mathbb{P}, \mathbb{Q}) \geq 1 - \sqrt{1 - \mathcal{S}_{\sqrt{}}^2(\mathbb{P}, \mathbb{Q})}.$$

PROOF. Introduce the variational distance $V(\mathbb{P}, \mathbb{Q})$ as the $f$-divergence associated with the convex function $f : u \mapsto \frac{1}{2}|u - 1|$. From Scheffé's theorem, we have $\mathcal{S}_\wedge(\mathbb{P}, \mathbb{Q}) = 1 - V(\mathbb{P}, \mathbb{Q})$ for any distributions $\mathbb{P}$ and $\mathbb{Q}$. Introduce the Hellinger distance $H$, which is defined as $H(\mathbb{P}, \mathbb{Q}) \geq 0$ and $1 - \frac{H^2(\mathbb{P},\mathbb{Q})}{2} = \mathcal{S}_\vee(\mathbb{P}, \mathbb{Q})$ for any probability distributions $\mathbb{P}$ and $\mathbb{Q}$. The variational and Hellinger distances are known (see, e.g., [46], Lemma 2.2) to be related by

$$V(\mathbb{P}, \mathbb{Q}) \leq \sqrt{1 - \left(1 - \frac{H^2(\mathbb{P}, \mathbb{Q})}{2}\right)^2},$$

hence the result. □

LEMMA 10.8. *For any distributions* $\mathbb{P}^{(1)}, \ldots, \mathbb{P}^{(k)}, \mathbb{Q}^{(1)}, \ldots, \mathbb{Q}^{(k)}$, *we have*

$$\mathcal{S}_\vee(\mathbb{P}^{(1)} \otimes \cdots \otimes \mathbb{P}^{(k)}, \mathbb{Q}^{(1)} \otimes \cdots \otimes \mathbb{Q}^{(k)})$$
$$= \mathcal{S}_\vee(\mathbb{P}^{(1)}, \mathbb{Q}^{(1)}) \times \cdots \times \mathcal{S}_\vee(\mathbb{P}^{(k)}, \mathbb{Q}^{(k)}).$$

PROOF. When it exists, the density of $\mathbb{P}^{(1)} \otimes \cdots \otimes \mathbb{P}^{(k)}$ w.r.t. $\mathbb{Q}^{(1)} \otimes \cdots \otimes \mathbb{Q}^{(k)}$ is the product of the densities of $\mathbb{P}^{(i)}$ w.r.t. $\mathbb{Q}^{(i)}$, $i = 1, \ldots, k$, hence the desired tensorization equality. □

From the last two lemmas, we obtain

$$\mathcal{S}_\wedge(Q_+^{\otimes k}, Q_-^{\otimes k}) \geq 1 - \sqrt{1 - S_\vee^{2k}(Q_+, Q_-)}.$$

Now we have

$$S_\vee^2(Q_+, Q_-) = [\sqrt{p_+ p_-} + \sqrt{(1 - p_+)(1 - p_-)}]^2$$
$$= 1 - [\sqrt{p_+(1 - p_-)} - \sqrt{(1 - p_+)p_-}]^2$$
$$= 1 - d_{\mathrm{II}}.$$

So we get

$$(10.18) \qquad \mathcal{S}_\wedge(Q_+^{\otimes k}, Q_-^{\otimes k}) \geq 1 - \sqrt{1 - (1 - d_{\mathrm{II}})^k} \geq 1 - \sqrt{k d_{\mathrm{II}}},$$

where the second inequality follows from the inequality $1 - x^k \leq k(1 - x)$ that holds for any $0 \leq x \leq 1$ and $k \geq 1$. This ends the proof of (10.17). □

By computing the second derivative of $u \mapsto \sqrt{1 - e^{-u}}$, we obtain that this function is concave. So for any $a \in [0; 1]$, the functions $x \mapsto 1 - \sqrt{1 - a^x}$ and $x \mapsto 1 - \sqrt{ax}$ are convex. The convexity of these functions and Lemmas 10.5 and 10.6 imply Lemma 10.4.

10.9. *Proofs of Theorems 8.3 and 8.4.* We consider a $(\tilde{m}, \tilde{w}, \tilde{d_{\mathrm{II}}})$-hypercube with

$$\tilde{m} = \lfloor \log_2 |\mathcal{G}| \rfloor,$$

$h_1 = -B$ and $h_2 = B$, and with $\tilde{w}$ and $\tilde{d_{\mathrm{II}}}$ to be taken in order to (almost) maximize the bound.

*Case $q = 1$.* Computations lead to

$$d_{\mathrm{I}} = |h_2 - h_1| \left[ \left| p_+ - \frac{1}{2} \right| \wedge \left| p_- - \frac{1}{2} \right| \right] = \frac{\sqrt{d_{\mathrm{II}}}}{2} |h_2 - h_1| = B\sqrt{d_{\mathrm{II}}}$$

so that, choosing $\tilde{w} = 1/\tilde{m}$, (8.14) gives

$$\sup_{P \in \mathcal{H}} \left\{ \mathbb{E} R(\hat{g}) - \min_g R(g) \right\} \geq B\sqrt{d_{\mathrm{II}}} (1 - \sqrt{n d_{\mathrm{II}}/\tilde{m}}).$$

Maximizing the lower bound w.r.t. $d_{\mathrm{II}}$, we choose $d_{\mathrm{II}} = \frac{\tilde{m}}{4n} \wedge 1$ and obtain the announced result.

*Case $1 < q \leq 1 + \sqrt{\frac{\tilde{m}}{4n} \wedge 1}$.* Tedious computations put in Appendix A.1 lead to: for any $p \in [0; 1]$,

(10.19) $$\phi(p) = p(1-p) \frac{|h_2 - h_1|^q}{[p^{1/(q-1)} + (1-p)^{1/(q-1)}]^{q-1}}$$

and

(10.20)
$$\phi''(p) = -\frac{q}{q-1} [p(1-p)]^{(2-q)/(q-1)}$$
$$\times \frac{|h_2 - h_1|^q}{[p^{1/(q-1)} + (1-p)^{1/(q-1)}]^{q+1}}.$$

From (8.11), for any $0 < \varepsilon \leq 1$, we get

$$d_{\mathrm{I}} \geq \frac{d_{\mathrm{II}}}{2} \int_{(1-\varepsilon)/2}^{(1+\varepsilon)/2} [t \wedge (1-t)] \left| \phi'' \left( \frac{1 - \sqrt{d_{\mathrm{II}}}}{2} + \sqrt{d_{\mathrm{II}}} t \right) \right| dt$$

$$\geq \frac{d_{\mathrm{II}}}{2} \frac{\varepsilon(2-\varepsilon)}{4} \inf_{u \in [(1-\varepsilon\sqrt{d_{\mathrm{II}}})/2; (1+\varepsilon\sqrt{d_{\mathrm{II}}})/2]} |\phi''(u)|$$

$$\geq \frac{\varepsilon(2-\varepsilon)}{8} d_{\mathrm{II}} \times \left| \phi'' \left( \frac{1 - \varepsilon\sqrt{d_{\mathrm{II}}}}{2} \right) \right|$$

$$\geq \frac{\varepsilon(2-\varepsilon)}{8} d_{\mathrm{II}} \times \frac{q}{q-1} \left[ \frac{1 - \varepsilon^2 d_{\mathrm{II}}}{4} \right]^{(2-q)/(q-1)} \frac{(2B)^q}{2^{q+1}[(1 + \varepsilon\sqrt{d_{\mathrm{II}}})/2]^{(q+1)/(q-1)}}$$

$$\geq \frac{\varepsilon(2-\varepsilon)}{8} d_{\mathrm{II}} \times \frac{4qB^q}{q-1}(1-\varepsilon\sqrt{d_{\mathrm{II}}})^{(2-q)/(q-1)}(1+\varepsilon\sqrt{d_{\mathrm{II}}})^{(1-2q)/(q-1)}$$

$$= (1-\varepsilon/2)qB^q\frac{\varepsilon d_{\mathrm{II}}}{q-1}(1-\varepsilon\sqrt{d_{\mathrm{II}}})^{(2-q)/(q-1)}(1+\varepsilon\sqrt{d_{\mathrm{II}}})^{(1-2q)/(q-1)}.$$

Let $K = (1-\varepsilon\sqrt{d_{\mathrm{II}}})^{(2-q)/(q-1)}(1+\varepsilon\sqrt{d_{\mathrm{II}}})^{(1-2q)/(q-1)}$. From (8.14), taking $\tilde{w} = 1/\tilde{m}$, we get

$$(10.21) \qquad \sup_{P\in\mathcal{H}}\left\{\mathbb{E}R(\hat{g})-\min_g R(g)\right\} \geq (1-\varepsilon/2)KqB^q\frac{\varepsilon d_{\mathrm{II}}}{q-1}(1-\sqrt{nd_{\mathrm{II}}/\tilde{m}}).$$

This leads us to choose $d_{\mathrm{II}} = \frac{\tilde{m}}{4n}\wedge 1$ and $\varepsilon = (q-1)\sqrt{\frac{n}{\tilde{m}}\vee\frac{1}{4}} \leq \frac{1}{2}$ and obtain

$$\mathbb{E}R(\hat{g})-\min_g R(g) \geq \frac{3qB^q}{8}K\left\{\left(\frac{1}{4}\sqrt{\frac{\tilde{m}}{n}}\right)\vee\left(1-\sqrt{\frac{n}{\tilde{m}}}\right)\right\}.$$

Since $1 < q \leq 2$ and $\varepsilon\sqrt{d_{\mathrm{II}}} = \frac{q-1}{2}$, we may check that $K \geq 0.29$ (to be compared with $\lim_{q\to 1}K = e^{-1}\approx 0.37$).

*Case $q > 1 + \sqrt{\frac{\tilde{m}}{4n}}$.* We take $\tilde{w} = \frac{1}{n+1}\wedge\frac{1}{\tilde{m}}$. From (8.4), (8.6) and (10.19), we get $d_{\mathrm{I}} = \psi_{1,0,-B,B}(1/2) = \phi_{-B,B}(1/2) = B^q$. From (8.15), we obtain

$$(10.22) \qquad \mathbb{E}R(\hat{g})-\min_{g\in\mathcal{G}} R(g) \geq \left(\frac{\lfloor\log_2|\mathcal{G}|\rfloor}{n+1}\wedge 1\right)B^q\left(1-\frac{1}{n+1}\wedge\frac{1}{\lfloor\log_2|\mathcal{G}|\rfloor}\right)^n$$

$$\geq e^{-1}B^q\left(\frac{\lfloor\log_2|\mathcal{G}|\rfloor}{n+1}\wedge 1\right),$$

where the last inequality uses $[1-1/(n+1)]^n \searrow e^{-1}$.

*Improvement when $1+\sqrt{\frac{\tilde{m}}{4n}}\wedge 1 < q < 2$.* From (10.21), by choosing $\varepsilon = 1/2$ and introducing $K' \triangleq (1-\sqrt{d_{\mathrm{II}}}/2)^{(2-q)/(q-1)}(1+\sqrt{d_{\mathrm{II}}}/2)^{(1-2q)/(q-1)}$, we obtain

$$\sup_{P\in\mathcal{H}}\left\{\mathbb{E}R(\hat{g})-\min_g R(g)\right\} \geq \frac{3qB^q}{8}K'\frac{d_{\mathrm{II}}}{q-1}(1-\sqrt{nd_{\mathrm{II}}/\tilde{m}}).$$

This leads us to choose $d_{\mathrm{II}} = \frac{4\tilde{m}}{9n}\wedge 1$. Since $\sqrt{\frac{\tilde{m}}{4n}}\wedge 1 < q-1$, we have $\sqrt{d_{\mathrm{II}}} \leq \frac{4}{3}(q-1)$, hence $K' \geq (1-\frac{2}{3}(q-1))^{(2-q)/(q-1)}(1+\frac{2}{3}(q-1))^{(1-2q)/(q-1)}$. For any $1 < q < 2$, this last quantity is greater than 0.2. So we have proved that for $1+\sqrt{\frac{\tilde{m}}{4n}}\wedge 1 < q < 2$,

$$(10.23) \qquad \mathbb{E}R(\hat{g})-\min_{g\in\mathcal{G}} R(g) \geq \frac{q}{90(q-1)}B^q\frac{\lfloor\log_2|\mathcal{G}|\rfloor}{n}.$$

Theorem 8.4 follows from (10.22) and (10.23).

### 10.10. *Proof of Theorem 8.6.*

10.10.1. *Proof of the first inequality of Theorem 8.6.* Let $\tilde{m} = \lfloor \log_2 |\mathcal{G}| \rfloor$. Contrary to other lower bounds obtained in this work, this learning setting requires asymmetrical hypercubes of distributions. Here we consider a constant $\tilde{m}$-dimensional hypercube of distributions with edge probability $\tilde{w}$ such that $p_+ \equiv p$, $p_- \equiv 0$, $h_1 \equiv +B$ and $h_2 \equiv 0$, where $\tilde{w}$, $p$ and $B$ are positive real parameters to be chosen according to the strategy described at the beginning of Section 8.3. To have $\mathbb{E}|Y|^s \leq A$, we need that $\tilde{m}\tilde{w}pB^s \leq A$. To ensure that a best prediction function has infinite norm bounded by $b$, from the computations at the beginning of Appendix A.1, we need that

$$B \leq \frac{p^{1/(q-1)} + (1-p)^{1/(q-1)}}{p^{1/(q-1)}} b.$$

This inequality is in particular satisfied for $B = Cp^{-1/(q-1)}$ for appropriate small constant $C$ depending on $b$ and $q$. From the definition of the edge discrepancy of type II, we have $d_{\mathrm{II}} = p$. In order to have the r.h.s. of (8.14) of order $mwd_{\mathrm{I}}$, we want to have $n\tilde{w}p \leq C < 1$. All the previous constraints lead us to take the parameters $\tilde{w}, p$ and $B$ such that

$$\begin{cases} B = Cp^{-1/(q-1)}, \\ \tilde{m}\tilde{w}pB^s = A, \\ n\tilde{w}p = 1/4. \end{cases}$$

Let $Q = \frac{\tilde{m}}{n} \wedge 1$. This leads to $p = CQ^{(q-1)/s}$, $B = CQ^{-1/s}$ and $\tilde{w} = C\tilde{m}^{-1}Q^{1-(q-1)/s}$ with $C$ small positive constants depending on $b$, $A$, $q$ and $s$. Now from the definition of the edge discrepancy of type I and (8.10), we have

$$\begin{aligned} d_{\mathrm{I}} &= \frac{p^2}{2} \int_0^1 [t \wedge (1-t)]|\phi_{0,B}''(tp)|\, dt \\ &\geq \frac{p^2}{2} \int_{1/4}^{3/4} \frac{1}{4} \min_{[p/4; 3p/4]} |\phi_{0,B}''(tp)|\, dt \\ &\geq Cp^2 p^{(2-q)/(q-1)} B^q \\ &= C, \end{aligned}$$

where the last inequality comes from (10.20). From (8.14), we get

$$\sup_{P \in \mathcal{P}} \left\{ \mathbb{E}R(\hat{g}) - \min_{g \in \mathcal{G}} R(g) \right\} \geq CQ^{1-(q-1)/s}.$$

10.10.2. *Proof of the second inequality of Theorem 8.6.* We still use $\tilde{m} = \lfloor \log_2 |\mathcal{G}| \rfloor$. We consider a $(\tilde{m}, \tilde{w}, \tilde{d_{\mathrm{II}}})$-hypercube with $h_1 \equiv -B$ and $h_2 \equiv +B$, where $\tilde{w}, \tilde{d_{\mathrm{II}}}$ and $B$ are positive real parameters to be chosen according to the strategy described at the beginning of Section 8.3. To have $\mathbb{E}|Y|^s \leq A$, we need that $\tilde{m}\tilde{w}B^s \leq A$. To ensure that a best prediction function has infinite norm bounded by $b$, from the computations at the beginning of Appendix A.1, we need that

$$(10.24) \qquad B \leq \frac{[1 + (\tilde{d_{\mathrm{II}}})^{1/2}]^{1/(q-1)} + [1 - (\tilde{d_{\mathrm{II}}})^{1/2}]^{1/(q-1)}}{[1 + (\tilde{d_{\mathrm{II}}})^{1/2}]^{1/(q-1)} - [1 - (\tilde{d_{\mathrm{II}}})^{1/2}]^{1/(q-1)}}\, b.$$

For fixed $q$ and $b$, this inequality essentially means that $B \leq C\tilde{d_{\mathrm{II}}}^{-1/2}$ since we intend to take $\tilde{d_{\mathrm{II}}}$ close to 0. In order to have the r.h.s. of (8.14) of order $mwd_{\mathrm{I}}$, we want to have $n\tilde{w}\tilde{d_{\mathrm{II}}} \leq 1/4$ where, once more, this last constant is arbitrarily taken. The previous constraints lead us to choose

$$\begin{cases} B = C\tilde{d_{\mathrm{II}}}^{-1/2}, \\ \tilde{m}\tilde{w}B^s = A, \\ n\tilde{w}\tilde{d_{\mathrm{II}}} = 1/4. \end{cases}$$

We still use $Q = \frac{\tilde{m}}{n} \wedge 1$. This leads to $\tilde{d_{\mathrm{II}}} = CQ^{2/(s+2)}$, $B = CQ^{-1/(s+2)}$ and $\tilde{w} = C\tilde{m}^{-1}Q^{s/(s+2)}$ with $C$ small positive constants depending on $b$, $A$, $q$ and $s$. Now from (10.20), we have $\phi''(t) \geq CB^q = CQ^{-q/(s+2)}$ for $t \in [p_-; p_+]$. Using (8.11) and (8.14), we obtain

$$\sup_{P \in \mathcal{P}} \left\{ \mathbb{E}R(\hat{g}) - \min_{g \in \mathcal{G}} R(g) \right\} \geq CQ^{1-q/(s+2)}.$$

## APPENDIX

**A.1. Computations of the second derivative of $\phi$ for the $L_q$-loss.** Let $h_1$ and $h_2$ be fixed. We start with the computation of $\phi$. For any $p \in [0; 1]$, the quantity $\varphi_p(y) = p|y - h_1|^q + (1-p)|y - h_2|^q$ is minimized when $y \in [h_1 \wedge h_2; h_1 \vee h_2]$ and $pq(y-h_1)^{q-1} = (1-p)q(h_2-y)^{q-1}$. Introducing $r = \frac{1}{q-1}$ and $D = p^r + (1-p)^r$, the minimizer can be written as $y = \frac{p^r h_1 + (1-p)^r h_2}{D}$ and the minimum is

$$\phi(p) = \left( p\frac{(1-p)^{rq}}{D^q} + (1-p)\frac{p^{rq}}{D^q} \right)|h_2 - h_1|^q$$

$$= p(1-p)\frac{|h_2 - h_1|^q}{D^{q-1}},$$

where we use the equality $rq = 1 + r$. We get

$$\frac{1}{|h_2 - h_1|^q}\phi'(p) = \frac{1-2p}{D^{q-1}} + p(1-p)(1-q)rD^{-q}[p^{r-1} - (1-p)^{r-1}]$$

$$= D^{-q}\{(1-2p)[p^r + (1-p)^r] - (1-p)p^r + p(1-p)^r]\}$$
$$= D^{-q}\{(1-p)^{r+1} - p^{r+1}\},$$

hence

$$\frac{1}{|h_2 - h_1|^q}\phi''(p) = -qrD^{-q-1}[p^{r-1} - (1-p)^{r-1}][(1-p)^{r+1} - p^{r+1}]$$

$$-qrD^{-q-1}[p^r - (1-p)^r]^2$$

$$= -qrD^{-q-1}p^{r-1}(1-p)^{r-1}$$

$$= -\frac{q}{q-1}\frac{[p(1-p)]^{(2-q)/(q-1)}}{[p^{1/(q-1)} + (1-p)^{1/(q-1)}]^{q+1}}.$$

**A.2. Expected risk bound from Hoeffding's inequality.** Let $\lambda' > 0$ and $\rho$ be a probability distribution on $\mathcal{G}$. Let $r(g)$ denote the empirical risk of a prediction function $g$, that is, $r(g) = \frac{1}{n}\sum_{i=1}^n L(Z_i, g)$. Hoeffding's inequality applied to the random variable $W = \mathbb{E}_{g\sim\rho}L(Z, g) - L(Z, g') \in [-(b-a); b-a]$ for a fixed $g'$ gives

$$\mathbb{E}_{Z\sim P}e^{\eta[W - \mathbb{E}W]} \le e^{\eta^2(b-a)^2/2}$$

for any $\eta > 0$. For $\eta = \lambda'/n$, this leads to

$$\mathbb{E}_{Z_1^n}e^{\lambda'[R(g') - \mathbb{E}_{g\sim\rho}R(g) - r(g') + \mathbb{E}_{g\sim\rho}r(g)]} \le e^{(\lambda')^2(b-a)^2/(2n)}.$$

Consider the Gibbs distribution $\hat{\rho} = \pi_{-\lambda'r}$. This distribution satisfies

$$\mathbb{E}_{g'\sim\hat{\rho}}r(g') + K(\hat{\rho}, \pi)/\lambda' \le \mathbb{E}_{g\sim\rho}r(g) + K(\rho, \pi)/\lambda'.$$

We have

$$\mathbb{E}_{Z_1^n}\mathbb{E}_{g'\sim\hat{\rho}}R(g') - \mathbb{E}_{g\sim\rho}R(g)$$

$$\le \mathbb{E}_{Z_1^n}\bigg\{\mathbb{E}_{g'\sim\hat{\rho}}[R(g') - \mathbb{E}_{g\sim\rho}R(g) - r(g') - \mathbb{E}_{g\sim\rho}r(g)]$$

$$+ \frac{K(\rho, \pi) - K(\hat{\rho}, \pi)}{\lambda'}\bigg\}$$

$$\le \frac{K(\rho, \pi)}{\lambda'} + \mathbb{E}_{Z_1^n}\frac{1}{\lambda'}\log\mathbb{E}_{g'\sim\pi}e^{\lambda'[R(g') - \mathbb{E}_{g\sim\rho}R(g) - r(g') - \mathbb{E}_{g\sim\rho}r(g)]}$$

$$\le \frac{K(\rho, \pi)}{\lambda'} + \frac{1}{\lambda'}\log\mathbb{E}_{g'\sim\pi}\mathbb{E}_{Z_1^n}e^{\lambda'[R(g') - \mathbb{E}_{g\sim\rho}R(g) - r(g') - \mathbb{E}_{g\sim\rho}r(g)]}$$

$$\le \frac{K(\rho, \pi)}{\lambda'} + \frac{\lambda'(b-a)^2}{2n}.$$

This proves that for any $\lambda > 0$, the generalization error of the algorithm which draws its prediction function according to the Gibbs distribution $\pi_{-\lambda\Sigma_n/2}$ satisfies

$$\mathbb{E}_{Z_1^n}\mathbb{E}_{g'\sim\pi_{-\lambda\Sigma_n/2}}R(g') \leq \min_{\rho\in\mathcal{M}}\left\{\mathbb{E}_{g\sim\rho}R(g) + 2\left[\frac{\lambda(b-a)^2}{8} + \frac{K(\rho,\pi)}{\lambda n}\right]\right\},$$

where we use the change of variable $\lambda = 2\lambda'/n$ in order to underline the difference with (6.4).

**Acknowledgments.**   I am very grateful to Nicolas Vayatis, Alexandre Tsybakov, Gilles Stoltz, Olivier Catoni and the reviewers for their helpful comments and suggestions.

Certis
Ecole des Ponts Paris Tech
Université Paris Est
6 Av Blaise Pascal—Cité Descartes
77455 Marne-la-Vallée
and
ENS/INRIA/CNRS, Willow project
Inria—Ecole Normale Supérieure
45 rue d'Ulm
75230 Paris
France
E-mail: audibert@certis.enpc.fr